\documentclass[11pt,letterpaper]{article}

\usepackage{amsmath,amssymb,amsfonts,amsthm}
\usepackage{tikz}
\usepackage{graphicx}
\usetikzlibrary{arrows.meta,positioning,calc}
\usepackage{float}
\usepackage{dsfont}
\usepackage{enumitem}
\usepackage{bbm}
\usepackage{xcolor}
\usepackage[ruled,vlined,linesnumbered]{algorithm2e}
\usepackage[margin=1in]{geometry}
\usepackage[colorlinks, citecolor = blue, linkcolor = magenta]{hyperref}
\usepackage[nameinlink]{cleveref}
\Crefname{ntn}{Notation}{Notations}
\Crefname{ques}{Question}{Questions}
\numberwithin{equation}{section}

\newtheorem{theorem}{Theorem}[section]
\newtheorem{lemma}[theorem]{Lemma}
\newtheorem{corollary}[theorem]{Corollary}

\theoremstyle{definition}

\newtheorem{proposition}[theorem]{Proposition}
\newtheorem{definition}[theorem]{Definition}

\newtheorem{notation}[theorem]{Notation}
\newtheorem{question}[theorem]{Question}

\Crefformat{question}{#2Question~#1#3}
\Crefformat{notation}{#2Notation~#1#3}

\renewcommand{\leq}{\leqslant}
\renewcommand{\le}{\leqslant}
\renewcommand{\geq}{\geqslant}

\newcommand{\cC}{\ensuremath{\mathcal{C}}}

\newcommand{\bC}{\ensuremath{\mathbb{C}}}

\newcommand{\bE}{\ensuremath{\mathbb{E}}}
\newcommand{\bF}{\ensuremath{\mathbb{F}}}

\newcommand{\bP}{\ensuremath{\mathbb{P}}}
\newcommand{\bR}{\ensuremath{\mathbb{R}}}

\newcommand{\bZ}{\ensuremath{\mathbb{Z}}}

\newcommand{\Ex}[1]{\bE \left[ #1 \right]}
\newcommand{\Exu}[2]{\underset{#1} \bE \left[ #2 \right] }

\newcommand{\Exs}[2]{\bE_{#1}\left[ #2 \right]}
\renewcommand{\Pr}[1]{\bP \left[ #1 \right]} % Probability

\newcommand{\Pru}[2]{\underset{ #1 }\bP \left[ #2 \right]}

\renewcommand{\d}{\mathrm{d}}
\newcommand{\ind}[1]{\mathds{1} \left[ #1 \right] }

\newcommand{\inp}[2]{\left\langle #1, #2 \right\rangle}

\title{Improved Almost-Periodicity Bounds and the Bogolyubov--Ruzsa Lemma}

\author{Yumou Fei\thanks{Massachusetts Institute of Technology, yumou415@mit.edu. Supported by NSF awards DMS-2022448 and CCF-2310818.}}

\date{}

\begin{document}

\maketitle

\begin{abstract}
The celebrated work of Sanders (\textit{Analysis \& PDE}, 2012), building on the almost-periodicity argument of Croot and Sisask (\textit{GAFA}, 2010), showed that if $A\subseteq\bF_2^n$ has density $\alpha>0$, then $2A-2A$ contains a linear subspace of codimension $O(\log(2/\alpha)^4)$. In this paper, we improve this upper bound on the codimension to $O(\log(2/\alpha)^3)$.

The main conceptual idea behind our proof is to maintain a distribution of almost-periods, instead of a set of almost-periods, in the argument of Croot and Sisask. This leads to a martingale structure in the subsequent convolution bootstrapping step, and the ``unbiased'' nature of the martingale saves a factor of $\log(2/\alpha)$ in the codimension bound.
\end{abstract}
\tableofcontents

\newpage 

\paragraph{AI statement.} The author worked extensively with ChatGPT Pro at all stages of this work. Despite being mathematically AI-assisted, all content in the paper is written with the author's own human understanding. We claim no credit for discovering the proof, and view our main contribution as presenting a way of understanding it. The author takes full responsibility for the paper's content.

\paragraph{Acknowledgement.} The author would like to thank Shuo Wang and Dor Minzer for stimulating discussions about the topic of this work.

\section{Introduction}

Let $G$ be a finite Abelian group, and let $\mu$ be a probability distribution on $G$. Starting from a point $x\in G$, we sample $a\sim\mu$ and move to $x-a$. The expected value of a function $f:G\to\bC$ after this single step is given by the convolution/averaging operator
\begin{equation}\label{eq:def_of_T_mu}
    (T_\mu f)(x):=\Exu{a\sim\mu}{f(x-a)}.
\end{equation}
The main (informal) question considered in this paper is:
\begin{question}\label{ques:main_1}
If $\mu$ is ``spread out'' on $G$, what can we say about the averaging effect of $T_\mu$?
\end{question}

To quantify how spread out $\mu$ is, we use the $L^q$-norms associated with the uniform probability measure on $G$.

\begin{notation}
Let $L^2(G)$ be the Hilbert space of all functions $f:G\to\bC$, equipped with the uniform probability measure. Throughout, the subscript $x\in G$ in a probability or expectation means that $x$ is sampled uniformly from $G$. Our inner product and $L^q$-norms are therefore
\[
    \langle f,g\rangle:=\Exu{x\in G}{f(x)\overline{g(x)}},
    \qquad
    \|f\|_q:=\left(\Exu{x\in G}{|f(x)|^q}\right)^{1/q}
    \quad (1\le q<\infty).
\]
\end{notation}

Let $h_\mu$ denote the probability density function\footnote{A probability density function over $G$ is a function $h:G\rightarrow [0,+\infty)$ such that $\Exs{x\in G}{h(x)}=1$.} of $\mu$, and we measure how spread out $\mu$ is by the $L^{2}$-norm $\|h_{\mu}\|_{2}$. Among all distributions on $G$, this norm ranges from $1$, attained only by the uniform distribution on $G$, to $|G|^{1/2}$, attained precisely by the point masses. We can now rewrite~\Cref{ques:main_1} as:

\begin{question}\label{ques:main_2}
What can an upper bound on $\|h_\mu\|_2$ tell us about the averaging effect of $T_\mu$?
\end{question}

The starting point of our exploration is a result in~\cite{croot2013arithmetic,schoen2016roth}. In order to present their result, we first introduce some convenient notation.

\begin{notation}
For any two functions $f,g\in L^{2}(G)$, we define the convolution $f*g\in L^{2}(G)$ as
\[
f*g(x):=\Exu{a\in G}{f(x-a)g(a)}.
\]
Therefore, for any distribution $\mu$, we have $T_{\mu}f=f*h_{\mu}$. We note that the convolution of two probability density functions is still a probability density function.
\end{notation}
\begin{notation}\label{ntn:operator_T_S}
Given a subset $S\subseteq G$, we let $\mathds{1}_{S}\in L^{2}(G)$ be the indicator function of $S$, defined by $\mathds{1}_{S}(x)=1$ for all $x\in S$ and $\mathds{1}_{S}(x)=0$ for all $x\in G\setminus S$. If $S$ is nonempty, let $h_{S}=(|G|/|S|)\cdot\mathds{1}_{S}$ be the probability density function of the uniform distribution on $S$. Accordingly, define $T_{S}:L^{2}(G)\rightarrow L^{2}(G)$ to be the operator that maps a function $f$ to $f*h_{S}$.
\end{notation}

\begin{notation}
Suppose $H$ is a subgroup of $G$. We say that a function $f:G\to\bC$ is $H$-invariant
if it is constant on each coset of $H$, or equivalently, if
$f(x+h)=f(x)$ for all $x\in G$ and $h\in H$. Note that the operator $T_{H}$ maps any function to an $H$-invariant function. 
\end{notation}

We can now state the relevant theorem from~\cite{schoen2016roth} as follows.\footnote{For simplicity, we discuss only the case $G=\bF_{2}^{n}$. The original statement of~\cite[Theorem 3.2]{schoen2016roth} treats the case $K_{1}=K_{2}$, but its proof readily adapts to separate parameters $K_{1}$ and $K_{2}$.}

\begin{theorem}[{\cite[Theorem 3.2]{schoen2016roth}}]\label{thm:SS16}
Let $n$ be a positive integer, and let $K_{1},K_{2}>1$ and $\varepsilon\in(0,1)$. Given any nonempty subsets $A,B\subseteq \bF_{2}^{n}$ with $\|h_{A}\|_{2}\leq 2^{K_{1}}$ and $\|h_{B}\|_{2}\leq 2^{K_{2}}$, for any subset $S\subseteq \bF_{2}^{n}$ there exists an $\bF_{2}$-linear subspace $V\subseteq \bF_{2}^{n}$ of codimension at most\footnote{The notation $O_{\varepsilon}(\cdot)$ hides constant factors that depend only on $\varepsilon$.} $O_{\varepsilon}(K_{1}K_{2}^{3})$ such that
\[
\left|\inp{T_{A}h_{B}-T_{V}(T_{A}h_{B})}{\mathds{1}_{S}}\right|\leq \varepsilon.
\]
\end{theorem}

Viewed from the perspective of~\Cref{ques:main_2}, this result studies the averaging effect of the operator $T_{A}$. Indeed,~\Cref{thm:SS16} says that $T_{A}h_{B}$ is close to a $V$-invariant function when tested against an arbitrary \emph{pre-determined} set $S$. This interpretation naturally raises several questions:
\begin{enumerate}[label=(\arabic*)]
\item Can we choose a single subspace $V$ of codimension $O_\varepsilon(K_{1}K_{2}^{3})$ that works simultaneously for every $S\subseteq\bF_2^n$? This is equivalent to asking for the bound
\[\|T_{A}h_{B}-T_{V}(T_{A}h_{B})\|_{1}\leq 2\varepsilon.\]
\item Can we further require that a single subspace $V$ works simultaneously for every nonempty $B\subseteq\bF_2^n$ satisfying $\|h_B\|_2\leq 2^{K_{2}}$? Since we are trying to characterize the averaging effect of the operator $T_{A}$, it would be ideal if the characterization depended not on the input function, but only on the operator itself. 
\item Can we replace $T_A$ by $T_\mu$ for an arbitrary distribution $\mu$ satisfying $\|h_\mu\|_2\leq 2^{K_{1}}$, and replace $h_B$ by an arbitrary complex-valued function $f\in L^2(G)$ satisfying $\|f\|_2\leq 2^{K_{2}}\|f\|_{1}$? Normalized indicator functions usually suffice for applications in additive combinatorics, but the operator viewpoint of~\Cref{ques:main_2} motivates asking how far these quantitative conclusions extend to general distributions and input functions.
\item Can we improve the upper bound on the codimension of $V$?
\end{enumerate}

Question (3) is the easiest to address: the original proof readily extends to this setting. It is also not hard to adapt the proof to meet the additional requirements in questions (1) and (2). We will show that, surprisingly, the natural adaptations that address questions (1)--(3) also create room for the quantitative improvement sought in question (4)! Specifically, we prove the following result.

\begin{theorem}\label{thm:main_F2}
Let $n$ be a positive integer, and let $K_{1},K_{2}>1$ and $\varepsilon\in (0,1)$. Given any probability distribution $\mu$ on $\bF_{2}^{n}$ with $\|h_{\mu}\|_{2}\leq 2^{K_{1}}$, there exists an $\bF_{2}$-linear subspace $V\subseteq \bF_{2}^{n}$ of codimension at most $O_{\varepsilon}(K_{1}K_{2}^{2})$ such that the following holds: for all $f\in L^{2}(\bF_{2}^{n})$ with $\|f\|_{2}\leq 2^{K_{2}}\|f\|_{1}$, we have
\begin{equation}\label{eq:main_F2}
\left\|T_{\mu}f-T_{V}(T_{\mu}f)\right\|_{1}\leq \varepsilon\|f\|_{1}.
\end{equation}
\end{theorem}

We stress that the subspace $V$ in~\Cref{thm:main_F2} does not depend on $f$, as required in question (2). The theorem also answers questions (1) and (3). As for question (4), it improves the bound on $\operatorname{codim} V$ in~\Cref{thm:SS16} from $O_{\varepsilon}(K_{1}K_{2}^{3})$ to $O_{\varepsilon}(K_{1}K_{2}^{2})$.

\subsection{Technical Overview}

For simplicity, we focus on the case $G=\bF_{2}^{n}$. Most of the results discussed here, both from the existing literature and from this paper, hold for general finite Abelian groups as well, with Bohr sets taking the place of subspaces (see~\Cref{def:bohr_set}). What is less convenient about general Abelian groups is that Bohr sets need not be subgroups.

The main idea behind the arguments of~\cite{croot2010probabilistic,sanders2012bogolyubov,croot2013arithmetic,schoen2016roth} leading to~\Cref{thm:SS16} is \emph{almost-periodicity}: informally, an element $a\in G$ is an almost-period of a function $f\in L^2(G)$ if the translation $T_{\{a\}}f$ is close to $f$ with respect to some metric.\footnote{Recall from~\Cref{ntn:operator_T_S} that $T_{\{a\}}$ is the translation operator $f\mapsto f(\cdot-a)$.} To prove~\Cref{eq:main_F2}, it suffices to show that every $a\in V$ is an $L^1$-almost-period of $T_\mu f$ with error at most $\varepsilon\|f\|_1$, since $T_V$ is the average of the operators $T_{\{a\}}$ over $a\in V$.

\begin{figure}[htbp]
    \centering
    \begin{tikzpicture}[
        box/.style={
            draw=black,
            solid,
            line width=.5pt,
            align=center,
            font=\small,
            text width=2.1cm,
            minimum height=11mm,
            inner xsep=2.5mm,
            inner ysep=2mm
        },
        perspective/.style={
            box,
            densely dotted,
            text width=3.8cm,
            minimum height=16mm
        },
        arrow/.style={
            -{Stealth[length=2mm,width=1.4mm]},
            line width=.55pt
        }
    ]
        % The two perspectives.
        \node[perspective] (a) at (0,1.1) {
            almost-periodicity\\
            from the physical-space\\
            perspective
        };
        \node[perspective] (b) at (0,-1.1) {
            almost-periodicity\\
            from the Fourier\\
            perspective
        };

        % The upper row.
        \node[box, text width=3cm, right=6mm of a] (one) {
            \Cref{subsec:empirical_approx,subsec:finding_almost_period,subsec:bootstrapping}
        };
        \node[box, fill=green!12, right=6mm of one] (two) {
            \Cref{subsec:martingale}
        };

        % The lower box lies directly below the martingale box.
        \node[box] (three) at (two.center |- b.center) {
            \Cref{subsec:Fourier_prelim}
        };

        % The final box lies halfway between the two rows.
        \node[box, anchor=west] (four)
            at ($(two.east)!0.5!(three.east)+(10mm,0)$) {
            \Cref{subsec:bohr_almost_periods}
        };

        % Dotted arrows from the two perspectives.
        \draw[arrow, densely dotted]
            (a.east) -- (one.west);
        \draw[arrow, densely dotted]
            (b.east) -- (three.west);

        % Solid arrows between the remaining boxes.
        \draw[arrow]
            (one.east) -- (two.west);
        \draw[arrow]
            (two.east) -- ([yshift=2.2mm]four.west);
        \draw[arrow]
            (three.east) -- ([yshift=-2.2mm]four.west);
    \end{tikzpicture}
    \caption{Outline of the argument. The main new ingredient,
    \Cref{subsec:martingale}, is highlighted in green.}
    \label{fig:proof_outline}
\end{figure}
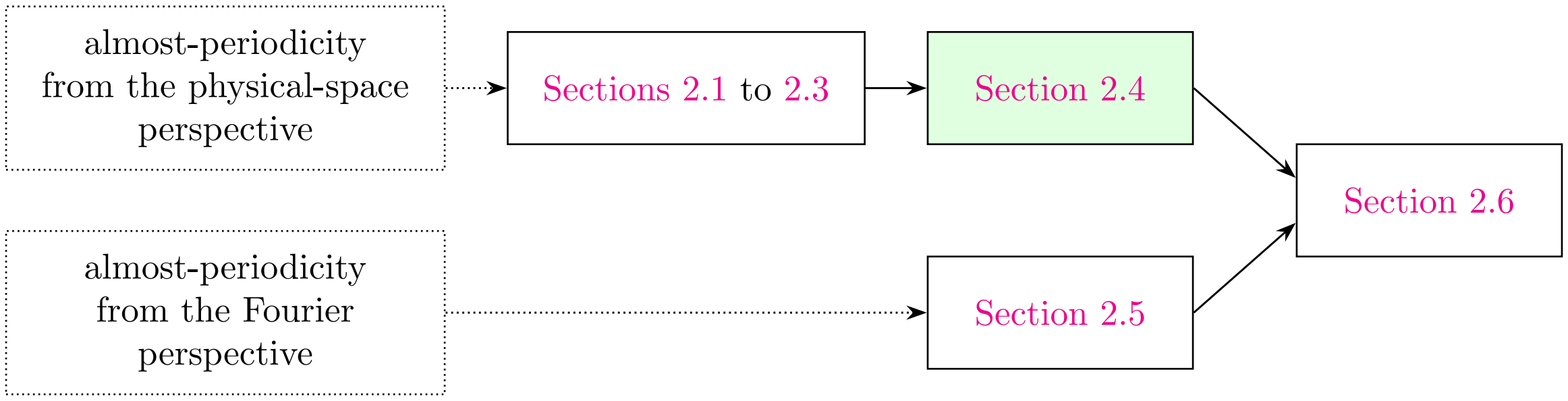

The probabilistic approach to almost-periodicity was pioneered by Croot and Sisask~\cite{croot2010probabilistic}, who showed how to obtain a large set of almost-periods of $T_\mu f$ when $\mu$ is uniform on a dense set $A\subseteq G$. Their method does not involve Fourier analysis: they use empirical approximation to find many approximators to $T_\mu f$, and whenever two of these approximators differ only by a translation, that translation is an almost-period of $T_\mu f$. By combining this result with Fourier-analytic arguments, Sanders~\cite{sanders2012bogolyubov} showed how to (roughly speaking) turn a large set of almost-periods into a subspace of almost-periods.

The main new conceptual message of this paper is that the empirical approximation method of~\cite{croot2010probabilistic} can be thought of as producing not a \emph{set} of almost-periods, but a \emph{distribution} of almost-periods for $T_{\mu}f$.\footnote{As with question (3) raised after~\Cref{thm:SS16}, passing from sets to distributions requires no essential changes to the subsequent argument, since the later steps extend naturally to distributions.} This distribution does not depend on $f$, while the almost-periodicity estimate holds in expectation for each $f$. As we explain in~\Cref{subsec:finding_almost_period}, this viewpoint naturally addresses questions (1) and (2) raised after~\Cref{thm:SS16}.

Perhaps surprisingly, because the distribution viewpoint preserves the ``unbiasedness'' inherent in the empirical approximation argument of~\cite{croot2010probabilistic}, it also reveals a martingale structure when we apply the bootstrapping argument of~\cite{sanders2012bogolyubov}. This yields a quantitative improvement in the upper bound on $\operatorname{codim} V$. The martingale step is the main new technical component of this work; it appears in~\Cref{subsec:martingale} and is highlighted in green in~\Cref{fig:proof_outline}.

In~\Cref{subsec:Fourier_prelim,subsec:bohr_almost_periods}, we then combine the resulting estimate, obtained entirely in ``physical space,'' with the Fourier-analytic argument of~\cite{sanders2012bogolyubov}. A small additional improvement is that we use a classical Bernstein-type inequality to control part of the $L^1$-error directly, instead of estimating it through an $L^2$-error that may be much larger. This improves the radius of the resulting Bohr set in general Abelian groups, but does not improve its rank or, in the case $G=\bF_2^n$, the codimension of the resulting subspace $V$.

\subsection{Applications to Additive Combinatorics}

Almost-periodicity has been a key ingredient in quantitative Bogolyubov--Ruzsa-type results in additive combinatorics~\cite{bogolyubov1939presque,chang2002polynomial,croot2010probabilistic,sanders2012bogolyubov,croot2013arithmetic,sanders2013structure,schoen2016roth}. It has also found important applications elsewhere in additive combinatorics; see, for example,~\cite{bourgain1990arithmetic,bloom2019logarithmic,kelley2023strong}.

We give two Bogolyubov--Ruzsa-type consequences of our almost-periodicity lemma (\Cref{lem:bohr_almost_period}). For a nonempty set $A\subseteq G$ of density $\alpha=|A|/|G|$, the first gives a Bohr set contained in $2A-2A$, while the second shows that $2A$ contains almost all of a translate of a Bohr set.

\begin{corollary}\label{cor:dense_bogolyubov}
Let $G$ be a finite Abelian group, and let $A\subseteq G$ be a nonempty set of density $\alpha=|A|/|G|$. Write $L=\log_2(2/\alpha)$. There exists a Bohr set $B=B(\Lambda,\rho)$ such that
\[
    B\subseteq2A-2A,
    \qquad
    |\Lambda|\leq CL^3,
    \qquad
    \rho\geq C^{-1}L^{-3},
\]
where $C>0$ is an absolute constant. In particular,
\(
    |B|\geq\exp\bigl(-C'L^3\log(2L)\bigr)|G|
\) for some absolute constant $C'>0$. If $G=\bF_2^n$, then $B$ can be taken to be a linear subspace of codimension at most $CL^3$.
\end{corollary}

\begin{corollary}\label{cor:dense_difference_bohr}
Let $G$ be a finite Abelian group, let $A\subseteq G$ be a nonempty set of density $\alpha=|A|/|G|$, and let $\varepsilon\in(0,1)$. Write $L=\log_2(2/\alpha)$. There exist a Bohr set $B=B(\Lambda,\rho)$ and an element $x\in 2A$ such that
\[
    |2A\cap(x+B)|\geq(1-\varepsilon)|B|,
    \qquad
    |\Lambda|\leq C_{\varepsilon}L^{3},
    \qquad
    \rho\geq C_{\varepsilon}^{-1}L^{-3},
\]
where $C_{\varepsilon}>0$ depends only on $\varepsilon$. In particular,
\(
|B|\geq
    \exp\bigl(-C_\varepsilon' L^3\log(2L)\bigr)|G|
\) for some $C'_\varepsilon>0$ depending only on $\varepsilon$. If $G=\bF_2^n$, then $x+B$ can be taken to be an affine subspace of codimension at most $C_{\varepsilon}L^{3}$.
\end{corollary}

The upper bounds on $|\Lambda|$ in~\Cref{cor:dense_bogolyubov,cor:dense_difference_bohr} improve the prior quartic-in-$L$ bounds of $O(L^4)$ and $O_{\varepsilon}(L^4)$, respectively, that follow from Sanders's work~\cite{sanders2012bogolyubov}; see also~\cite[Section~9]{schoen2016roth}.

\section{Almost-Periodicity}

In this section, we prove the main almost-periodicity lemma of this paper~(\Cref{lem:bohr_almost_period}). We then derive~\Cref{thm:main_F2} and~\Cref{cor:dense_bogolyubov,cor:dense_difference_bohr} from it in~\Cref{sec:applications}.

Throughout this section, let $G$ be a finite Abelian group. We occasionally use the case $G=\bF_2^n$ to guide the intuition, but all formal definitions and statements apply to general finite Abelian groups.

We assume the reader is familiar with Fourier analysis on finite Abelian groups. Apart from this, the exposition is largely self-contained. The martingale estimates used in~\Cref{subsec:martingale} are proved in~\Cref{subsec:martingale_increment}, and the Fourier-analytic tools introduced in~\Cref{subsec:Fourier_prelim} are proved in~\Cref{sec:Fourier_tools}.

\subsection{Empirical Approximation}\label{subsec:empirical_approx}

Given a probability distribution $\mu$ on $G$, a standard way to obtain a coarse empirical approximation is to draw $m$ independent samples $a_1,\dots,a_m\sim\mu$ and form the distribution $\nu$ with density
\[
    h_\nu(x)=\frac{|G|}{m}\sum_{i=1}^m\ind{x=a_i}.
\]
It is not hard to see that if both $\|h_{\mu}\|_{2}$ and $m$ are much smaller than $|G|^{1/4}$, then such an empirical distribution $\nu$ is far from $\mu$ in total variation. The first key observation of~\cite{croot2010probabilistic}, formalized in~\Cref{prop:T_nu_approximates_T_mu}, is that $T_\nu f$ can nevertheless approximate $T_\mu f$ well.

\begin{notation}
Given a positive integer $m$ and a tuple $Y=(y_1,\dots,y_m)\in G^m$, define the averaging operator $T_Y:L^2(G)\to L^2(G)$ by
\[
    T_Y=\frac{1}{m}\sum_{i=1}^m T_{\{y_i\}}.
\]
\end{notation}

\begin{proposition}\label{prop:T_nu_approximates_T_mu}
Let $m$ be a positive integer, and let $\mu$ be a probability distribution on $G$. For every $f\in L^2(G)$ and every even integer $q\geq 2$, we have\footnote{Here $\mu^{\otimes m}$ denotes the product distribution on $G^m$, whose coordinates are independent and each distributed according to $\mu$.}
\[
    \Exu{Y\sim\mu^{\otimes m}}{
        \big\|T_\mu f-T_Yf\big\|_{q}^{q}
    }
    \leq
    \left(
        \sqrt{\frac{2q}{m}}\,\|f\|_{q}
    \right)^{q}.
\]
\end{proposition}

\begin{proof}
Fix an element $x\in G$. The random variables $f(x-y_1),\dots,f(x-y_m)$ are independent and each has expectation $(T_\mu f)(x)$. Applying~\Cref{prop:Marcinkiewicz-Zygmund} and then H\"older's inequality gives
\begin{align*}
    \Exu{Y\sim\mu^{\otimes m}}{
        \Big|(T_Yf)(x)-(T_\mu f)(x)\Big|^{q}
    }
    &\leq
    \frac{(2q)^{q/2}}{m^{q}}\cdot
    \Exu{Y\sim\mu^{\otimes m}}{
        \left(\sum_{i=1}^m \Big|f(x-y_i)\Big|^2\right)^{q/2}
    }\\
    &\leq
    \frac{(2q)^{q/2}}{m^{q/2+1}}\cdot 
    \Exu{Y\sim\mu^{\otimes m}}{
        \sum_{i=1}^m \Big|f(x-y_i)\Big|^{q}
    }\\
    &=
    \left(\frac{2q}{m}\right)^{q/2}
    \big(|f|^{q}*h_{\mu}\big)(x).
\end{align*}
Taking expectation over $x\in G$ and using
\[
    \Exu{x\in G}{\big(|f|^{q}*h_{\mu}\big)(x)}
    =
    \Exu{x\in G}{|f(x)|^{q}}
    =
    \|f\|_{q}^{q}
\]
yields the desired conclusion.
\end{proof}

\subsection{Finding Almost-Periods}\label{subsec:finding_almost_period}

Applying Markov's inequality to~\Cref{prop:T_nu_approximates_T_mu} yields a large set of tuples $Y$ for which $T_Yf$ is close to $T_\mu f$. The next crucial observation of~\cite{croot2010probabilistic} is that, if $\mu$ is sufficiently spread out, then this set would be so large that there would be two distinct such tuples $Y,Y'$ differing by a common translation in every coordinate. More precisely, this means there is a common $a\in G$ such that $y_{i}-y_{i}'=a$ for all $i\in [m]$, and hence $T_Y=T_{\{a\}}T_{Y'}$. Since both $T_{Y}f$ and $T_{Y'}f$ are close to $T_{\mu}f$, the triangle inequality shows that $T_\mu f$ is close to $T_{\{a\}}(T_{\mu}f)$, and hence $a$ is an almost-period of $T_{\mu}f$.

\subsubsection{A Set of Almost-Periods}

\label{subsubsec:set_almost_periods}

In~\Cref{subsubsec:set_almost_periods}, we briefly explain how the argument of~\cite{croot2010probabilistic} produces a \emph{set} of almost-periods; this part will be superseded by~\Cref{subsubsec:distribution_almost_periods} and will not be used later in the paper.

For simplicity, within~\Cref{subsubsec:set_almost_periods} we suppose $\mu$ is uniform on a nonempty set $A\subseteq G$, and write $\alpha=|A|/|G|$. Fix $f\in L^2(G)$, a positive integer $m$, and an even integer $q\geq 2$. By~\Cref{prop:T_nu_approximates_T_mu} and Markov's inequality, the set
\[
    \mathcal{Y}
    :=
    \left\{
        Y\in A^m\,\middle|\,
        \big\|T_Af-T_Yf\big\|_q
        \leq
        2\sqrt{2q/m}\,\|f\|_q
    \right\}
\]
has size at least $|A|^m/2=\alpha^{m}|G|^{m}/2$.

Partition $G^m$ into the $|G|^{m-1}$ cosets of its diagonal subgroup $\{(a,\dots,a)\mid a\in G\}$. By the pigeonhole principle, some coset contains at least $\alpha^m|G|/2$ elements of $\mathcal{Y}$. Fix one such element $Y_0$, and let
\[
    S:=\{a\in G:Y_0+(a,\dots,a)\in\mathcal{Y}\}.
\]
Then $0\in S$ and $|S|\geq\alpha^m|G|/2$. For any $a\in S$, the corresponding tuple $Y=Y_0+(a,\dots,a)$ satisfies $T_Y=T_{\{a\}}T_{Y_{0}}$. We then have
\begin{align*}
\big\|T_{A} f-T_{\{a\}}(T_{A} f)\big\|_{q}
    &\leq
    \big\|T_{A} f-T_Yf\big\|_{q}
    +
    \big\|T_Yf-T_{\{a\}}(T_A f)\big\|_{q}\tag{by the triangle inequality}\\
    &=
    \big\|T_A f-T_Yf\big\|_{q}
    +
    \big\|T_A f-T_{Y_{0}}f\big\|_{q}\tag{translation invariance of the $L^{q}$-norm}\\
    &\leq 2\sqrt{2q/m}\,\|f\|_{q}+2\sqrt{2q/m}\,\|f\|_{q}\tag{since both $Y$ and $Y_{0}$ belong to $\mathcal{Y}$}\\
    &=4\sqrt{2q/m}\,\|f\|_{q}.
\end{align*}
Thus every element of $S$ is an almost-period of $T_Af$. Averaging over $S$, we obtain
\[
    \big\|T_Af-T_S(T_Af)\big\|_q
    \leq
    4\sqrt{\frac{2q}{m}}\,\|f\|_q.
\]

\subsubsection{A Distribution of Almost-Periods}\label{subsubsec:distribution_almost_periods}

As explained in~\Cref{subsubsec:set_almost_periods}, the argument of~\cite{croot2010probabilistic} finds a set $S$ of almost-periods for the function $T_{\mu}f$ (at least when $\mu$ is uniform over a set). However, the set $S$ of almost-periods depends on $f$, while, as mentioned in the introduction, we want the distribution of almost-periods to depend only on $\mu$. With this in mind, we deviate from~\cite{croot2010probabilistic} at this step. We note that the dependence on $f$ is introduced by selecting via Markov's inequality the tuples $Y$ (sampled from $\mu^{\otimes m}$) for which $T_Yf$ is close to $T_\mu f$. Therefore, refraining from applying Markov's inequality, we work directly with the full distribution $\mu^{\otimes m}$ and consider the distribution of translations defined below.
\begin{definition}\label{def:D_mu_m}
Given a positive integer $m$ and a probability distribution $\mu$ on $G$, define $D(\mu,m)$ to be the distribution on $G$ generated by the following procedure.
\begin{enumerate}
    \item Sample $Y=(y_1,\dots,y_m)$ from $\mu^{\otimes m}$.
    \item With $Y$ fixed, sample $Y'=(y_1',\dots,y_m')$ from $\mu^{\otimes m}$ conditioned on
    \[
        y_1-y_1'=y_2-y_2'=\cdots=y_m-y_m'.
    \]
    \item Output $a=y_1-y_1'$.
\end{enumerate}

\end{definition}

\begin{proposition}\label{prop:symmetry_D_mu_m}
In the notation and setting of~\Cref{def:D_mu_m}, the joint distribution of $(Y,Y')$ is invariant under interchanging $Y$ and $Y'$. In particular, the marginal distribution of $Y'$ is $\mu^{\otimes m}$, the same as that of $Y$.
\end{proposition}
\begin{proof}
We partition the group $G^{m}$ into cosets of the diagonal subgroup $\{(a,a,\dots,a)\mid a\in G\}$. The procedure in~\Cref{def:D_mu_m} first samples $Y$ according to $\mu^{\otimes m}$, and then resamples $Y'$ within the coset of $Y$. Therefore, conditioned on the coset, $Y$ and $Y'$ are independent samples from the same conditional distribution.   
\end{proof}

We can now prove the following almost-periodicity result for $T_{\mu}f$.

\begin{proposition}\label{prop:almost_period}
Let $m$ be a positive integer, and let $\mu$ be a probability distribution on $G$. For every $f\in L^2(G)$ and every even integer $q\geq 2$, we have
\[
    \Exu{a\sim D(\mu,m)}{
        \big\|T_\mu f-T_{\{a\}}(T_\mu f)\big\|_{q}
    }
    \leq
    \sqrt{\frac{8q}{m}}\,\|f\|_{q}.
\]
\end{proposition}

\begin{proof}
Let $Y,Y'$, and $a$ be sampled as in the definition of $D(\mu,m)$. Since $y_i=y_i'+a$ for every $i$, we have
\[
    T_Yf=T_{\{a\}}(T_{Y'}f).
\]
The triangle inequality and translation invariance of the $L^{q}$-norm give
\[
    \big\|T_\mu f-T_{\{a\}}(T_\mu f)\big\|_{q}
    \leq
    \big\|T_\mu f-T_Yf\big\|_{q}
    +
    \big\|T_Yf-T_{\{a\}}(T_\mu f)\big\|_{q}=
    \big\|T_\mu f-T_Yf\big\|_{q}
    +
    \big\|T_\mu f-T_{Y'}f\big\|_{q}.
\]
Taking expectations and using the common marginal distribution of $Y$ and $Y'$, we obtain
\begin{align*}
    \Exu{a\sim D(\mu,m)}{
        \big\|T_\mu f-T_{\{a\}}(T_\mu f)\big\|_{q}
    }
    &\leq
    2\cdot \Exu{Y\sim\mu^{\otimes m}}{
        \big\|T_\mu f-T_Yf\big\|_{q}
    }\\
    &\leq
    2\left(
        \Exu{Y\sim\mu^{\otimes m}}{
            \big\|T_\mu f-T_Yf\big\|_{q}^{q}
        }
    \right)^{1/q}\\
    &\leq
    \sqrt{\frac{8q}{m}}\,\|f\|_{q},
\end{align*}
where the second inequality is H\"older's inequality and the last follows from~\Cref{prop:T_nu_approximates_T_mu}.
\end{proof}
\begin{corollary}\label{cor:almost_period}
In the setting and notation of~\Cref{prop:almost_period}, we have
\begin{equation}\label{eq:cor_almost_period}
\left\|T_{\mu}f-T_{D(\mu,m)}(T_{\mu}f)\right\|_{q}\leq \sqrt{\frac{8q}{m}}\,\|f\|_{q}.
\end{equation}
\end{corollary}
\begin{proof}
The conclusion follows directly from~\Cref{prop:almost_period} by the convexity of the $L^{q}$-norm.
\end{proof}

\subsection{Bootstrapping into Subspace-Invariance}\label{subsec:bootstrapping}

The conclusion of~\Cref{cor:almost_period} already has a form somewhat similar to the desired bound in~\Cref{eq:main_F2} of~\Cref{thm:main_F2}. There are still two main issues:
\begin{itemize}
    \item We would like to replace the operator $T_{D(\mu,m)}$ in~\Cref{eq:cor_almost_period} by $T_V$ for some linear subspace $V\subseteq\bF_2^n$.
    \item We need to bound the codimension of $V$.
\end{itemize}
For the second item, we note that bounding the codimension of $V$ is equivalent to bounding $\|h_V\|_2$, since $\|h_V\|_2=2^{\operatorname{codim}(V)/2}$. Before addressing the first item, namely replacing $T_{D(\mu,m)}$ by $T_V$, we note that an upper bound on $\|h_\mu\|_2$ gives an upper bound on $\|h_{D(\mu,m)}\|_2$.

\begin{lemma}\label{lem:D_mu_m_density}
Given a positive integer $m$ and a probability distribution $\mu$ on $G$, we have
\[
    \left\|h_{D(\mu,m)}\right\|_2
    \leq
    \left\|h_\mu\right\|_2^m.
\]
\end{lemma}
\begin{proof}
Let $\cC$ be the collection of the $|G|^{m-1}$ cosets of the diagonal subgroup of $G^m$. For $Y=(y_1,\dots,y_m)\in G^m$, write
\[
    p(Y):=\frac{1}{|G|^m}\prod_{i=1}^m h_\mu(y_i)
\]
for the probability of sampling $Y$ from $\mu^{\otimes m}$.

For each $Y$ with $p(Y)>0$, let $\nu_Y$ be the conditional distribution of the output $a$ in~\Cref{def:D_mu_m}, given this choice of $Y$. If $C\in \cC$ is the coset that $Y$ belongs to, then $Y'$ is sampled from $\mu^{\otimes m}$ conditioned on belonging to $C$. As $Y'$ ranges over $C$, the output $a=y_1-y_1'$ ranges bijectively over $G$. Therefore,
\[
    \left\|h_{\nu_Y}\right\|_2
    =\frac{\left(\Exs{Z\in C}{p(Z)^{2}}\right)^{1/2}}{\Exs{Z\in C}{p(Z)}}=
    \frac{|G|^{1/2}}{\mu^{\otimes m}(C)}
    \left(\sum_{Z\in C}p(Z)^2\right)^{1/2}.
\]

Since $D(\mu,m)$ is the average of these conditional distributions over $Y\sim\mu^{\otimes m}$, the convexity of the $L^2$-norm gives
\begin{align*}
    \left\|h_{D(\mu,m)}\right\|_2
    &\leq
    \Exu{Y\sim\mu^{\otimes m}}{\left\|h_{\nu_Y}\right\|_2}=
    |G|^{1/2}\sum_{C\in\cC}
        \left(\sum_{Z\in C}p(Z)^2\right)^{1/2}\\
    &\leq
    \left(
        |G|^m\sum_{Y\in G^m}p(Y)^2
    \right)^{1/2}=
    \|h_\mu\|_2^m,
\end{align*}
where the second inequality is Cauchy--Schwarz over the $|G|^{m-1}$ cosets, and the last equality follows from the product structure of $p$.
\end{proof}

Therefore, if we can turn $T_{D(\mu,m)}$ into $T_V$ for some $V$, while also using the upper bound on $\left\|h_{D(\mu,m)}\right\|_2$ to control $\|h_V\|_2$, we will address both issues mentioned above.

The key idea, following~\cite{sanders2012bogolyubov}, is to consider a high convolution power of $D(\mu,m)$. Under repeated convolution, Fourier coefficients of small magnitude decay rapidly. Chang's lemma~\cite{chang2002polynomial} then controls the dimension of the span of the remaining large coefficients, leading to a distribution that is close to $V$-invariant for some subspace $V$ of small codimension.

\begin{notation}
Given a probability distribution $\nu$ on $G$ and a positive integer $k$, we write $\nu^{*k}$ for the $k$-fold convolution of $\nu$ with itself, or equivalently, the distribution of $a_1+\cdots+a_k$, where $a_1,\dots,a_k$ are independent samples from $\nu$. We write $T_\nu^k$ for the $k$-fold composition of $T_\nu$ with itself. Thus,
\[
    T_\nu^k=T_{\nu^{*k}},
    \qquad
    h_{\nu^{*k}}
    =
    \underbrace{h_\nu*\cdots*h_\nu}_{k\text{ factors}}
    =:
    (h_\nu)^{*k}.
\]
\end{notation}

We can now iterate~\Cref{cor:almost_period} to obtain the following bound.

\begin{proposition}\label{prop:k-fold_almost_period}
Let $m,k$ be positive integers, and let $\mu$ be a probability distribution on $G$. For every $f\in L^{2}(G)$ and every even integer $q\geq 2$, we have
\begin{equation}\label{eq:k-fold_almost_period}
\left\|T_{\mu}f-T_{D(\mu,m)}^{k}(T_{\mu}f)\right\|_{q}\leq k\sqrt{\frac{8q}{m}}\,\|f\|_{q}.
\end{equation}
\end{proposition}
\begin{proof}
This follows immediately from~\Cref{cor:almost_period} and the following general observation: for any probability distribution $\nu$ on $G$ and any $g\in L^2(G)$, we have
\begin{equation}\label{eq:telescope_inefficient}
    \left\|g-T_\nu^k g\right\|_{q}
    \leq
    k\left\|g-T_\nu g\right\|_{q}.
\end{equation}
To see this, apply the triangle inequality to the telescoping identity
\[
    g-T_\nu^k g
    =
    \sum_{j=0}^{k-1}
    T_\nu^j(g-T_\nu g)
\]
and use the fact that $T_\nu$ is a contraction in the $L^{q}$-norm.
\end{proof}

Starting from~\Cref{prop:k-fold_almost_period}, one can continue to follow the argument in~\cite{sanders2012bogolyubov,croot2013arithmetic} to obtain a version of~\Cref{thm:main_F2} with the weaker upper bound $\operatorname{codim} V\leq O_{\varepsilon}(K_{1}K_{2}^{3})$. However, we choose not to proceed at this point, because it turns out that the factor of $k$ on the right-hand side of~\Cref{eq:k-fold_almost_period} can be improved to an $O(\sqrt{k})$ factor! We will first establish this improvement in the next subsection, and then return to the argument of~\cite{sanders2012bogolyubov,croot2013arithmetic} after that.

\subsection{The Martingale Structure}\label{subsec:martingale}

As before, let $q\geq 2$ be an even integer. Consider a simple scenario where we have $k$ independent random variables $X_{1},\dots,X_{k}$, each with a $q$-norm $\Ex{|X_{i}|^{q}}^{1/q}$ at most 1. By the triangle inequality, the $q$-norm of $X_{1}+\dots+X_{k}$ is at most $k$. However, if $X_{1},\dots,X_{k}$ are all \emph{unbiased} variables, i.e. have mean 0, then the $q$-norm of their sum is actually at most $O(\sqrt{kq})$. This follows from~\Cref{prop:Marcinkiewicz-Zygmund}, which we have already used in the proof of~\Cref{prop:T_nu_approximates_T_mu}. Intuitively, when the variables $X_i$ are unbiased, their sum accumulates ``fluctuations'' rather than ``systematic biases,'' and fluctuations can cancel with each other. That is why, for fixed $q$, accumulating $k$ such fluctuations gives only an $O(\sqrt{k})$ dependence on $k$.

Looking back into the proof of~\Cref{prop:k-fold_almost_period}, one sees that the argument for~\Cref{eq:telescope_inefficient} uses the triangle inequality and does not account for any cancellation. Can we analyze the $k$ successive applications of $T_{D(\mu,m)}$ as an accumulation of mean-zero fluctuations?

If we proceed as in~\cite{croot2010probabilistic}, finding a \emph{set} $S$ of almost-periods as discussed in~\Cref{subsubsec:set_almost_periods}, then we first use Markov's inequality to select tuples $Y$ for which $T_Yf$ is close to $T_\mu f$. Before this selection, the empirical approximation is unbiased (see the proof of~\Cref{prop:T_nu_approximates_T_mu}):
\[
    \Exs{Y\sim\mu^{\otimes m}}{T_Yf}=T_\mu f.
\]
However, restricting to the selected tuples need not preserve this identity. Thus, the resulting approximation errors need not remain unbiased. In contrast, during our construction of the \emph{distribution} $D(\mu,m)$ of almost-periods, we did not discard any tuples from the support of $\mu^{\otimes m}$. The ``unbiasedness'' that is inherent in $\mu^{\otimes m}$ is therefore preserved, and this turns out to be the key for us to replace the factor of $k$ in~\Cref{eq:k-fold_almost_period} by an $O(\sqrt{k})$ factor!

\begin{lemma}\label{lem:martingale_almost_period}
Let $m,k$ be positive integers, and let $\mu$ be a probability distribution on $G$. For every $f\in L^{2}(G)$ and every real number $q\geq 2$, we have
\[
\left\|T_{\mu}f-T_{D(\mu,m)}^{k}(T_{\mu}f)\right\|_{q}\leq 4\sqrt{\frac{qk}{m}}\,\|f\|_{q}.
\]
\end{lemma}

\begin{proof}
Consider the following stochastic process:
\begin{enumerate}
    \item Sample $x_0\in G$ uniformly at random.
    \item Repeat the following for $j=1,2,\dots,k$:
    \begin{enumerate}
        \item Independently of the preceding history, sample $Y_{(j)}=(y_{j,1},\dots,y_{j,m})$ from $\mu^{\otimes m}$.
        \item Let
        \[
            Z_{(j)}
            =
            (z_{j,1},\dots,z_{j,m})
            =
            (x_{j-1}-y_{j,1},\dots,x_{j-1}-y_{j,m})
            \in G^m.
        \]
        \item With $Y_{(j)}$ fixed, sample $Y_{(j)}'=(y_{j,1}',\dots,y_{j,m}')$ from $\mu^{\otimes m}$ conditioned on
        \[
            y_{j,1}-y_{j,1}'
            =
            y_{j,2}-y_{j,2}'
            =
            \cdots
            =
            y_{j,m}-y_{j,m}'.
        \]
        \item Let $x_j=x_{j-1}-y_{j,1}+y_{j,1}'$.
    \end{enumerate}
\end{enumerate}
As illustrated in~\Cref{fig:resampling_process}, this process alternates between the spaces $G$ and $G^m$.

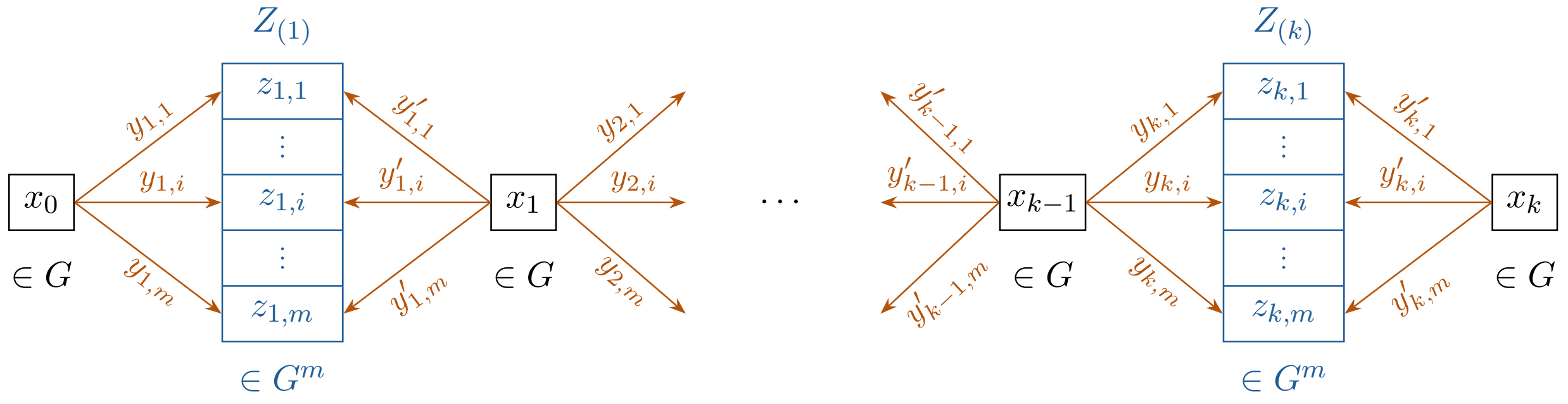
\begin{figure}[H]
    \centering
    \definecolor{processZ}{HTML}{1F5D99}
    \definecolor{processY}{HTML}{B45309}
    \resizebox{\linewidth}{!}{%
\begin{tikzpicture}[
    x=1cm,y=1cm,
    state/.style={draw=black,text=black,line width=.5pt,minimum width=7mm,
                  minimum height=6mm,inner sep=2pt},
    entry/.style={text=processZ,minimum width=13mm,minimum height=6mm,
                  inner sep=0pt},
    arrow/.style={draw=processY,-{Stealth[length=1.8mm,width=1.3mm]},
                  line width=.5pt},
    edge label/.style={text=processY,font=\small,fill=none,draw=none,inner sep=2pt}
]
    % The states and the omitted intermediate steps.
    \node[state] (x0)  at (0,0)    {$x_0$};
    \node[state] (x1)  at (5.2,0)  {$x_1$};
    \node         at (8.0,0)      {$\cdots$};
    \node[state] (xkm) at (10.8,0) {$x_{k-1}$};
    \node[state] (xk)  at (16.0,0) {$x_k$};
    \foreach \name in {x0,x1,xkm,xk}
        \node[anchor=north] at ([yshift=-2mm]\name.south)
            {$\in G$};

    % Draw the first and last tuple columns. The same pattern applies
    % to every intermediate step j.
    \foreach \j/\xc/\leftstate/\rightstate in
        {1/2.6/x0/x1,k/13.4/xkm/xk} {
        \node[text=processZ] at (\xc,1.9) {$Z_{(\j)}$};
        \node[text=processZ] at (\xc,-1.9) {$\in G^m$};

        % A vertical grid, with omitted coordinates represented by dots.
        \draw[draw=processZ,line width=.5pt]
            ({\xc-.65},-1.5) rectangle ({\xc+.65},1.5);
        \foreach \yy in {-.9,-.3,.3,.9}
            \draw[draw=processZ,line width=.5pt]
                ({\xc-.65},\yy) -- ({\xc+.65},\yy);
        % Center the visible dots geometrically in each cell.
        \foreach \yc in {.6,-.6}
            \foreach \dy in {-.10,0,.10}
                \fill[processZ] (\xc,{\yc+\dy}) circle[radius=.45pt];

        % Each arrow points from a state to a tuple entry.
        \foreach \i/\yy/\side in {1/1.2/above,i/0/above,m/-1.2/below} {
            \node[entry] (z-\j-\i) at (\xc,\yy)
                {$z_{\j,\i}$};
            \draw[arrow] (\leftstate.east) --
                node[edge label,pos=.60,sloped,\side]
                    {$y_{\j,\i}$}
                (z-\j-\i.west);
            \draw[arrow] (\rightstate.west) --
                node[edge label,pos=.60,sloped,\side]
                    {$y_{\j,\i}'$}
                (z-\j-\i.east);
        }
    }

    % Arrows to the omitted entries of Z_(2) and Z_(k-1).
    \foreach \i/\yy/\side in {1/1.2/above,i/0/above,m/-1.2/below} {
        \draw[arrow] (x1.east) --
            node[edge label,pos=.60,sloped,\side]
                {$y_{2,\i}$}
            (6.95,\yy);
        \draw[arrow] (xkm.west) --
            node[edge label,pos=.60,sloped,\side]
                {$y_{k-1,\i}'$}
            (9.05,\yy);
    }
\end{tikzpicture}
    }
    \caption{The stochastic process. The arrows encode the relations
    $z_{j,i}=x_{j-1}-y_{j,i}=x_j-y_{j,i}'$.
    }
    \label{fig:resampling_process}
\end{figure}

We make the following important observations about this process:

\begin{enumerate}[label=(\arabic*)]
\item Comparing with~\Cref{def:D_mu_m}, we see that, for each $j\in[k]$, conditioned on $x_0,\dots,x_{j-1}$, the distribution of $x_{j-1}-x_j$ is $D(\mu,m)$. In particular, every $x_j$ is uniformly distributed on $G$. Induction on $j$ also gives
    \[
        \Ex{(T_\mu f)(x_j)\,\big|\,x_0}
        =
        \left(T_{D(\mu,m)}^j(T_\mu f)\right)(x_0)
        \qquad\text{for all }j\in[k].
    \]
Consequently,
\begin{align}
\left\|T_{\mu}f-T_{D(\mu,m)}^{k}(T_{\mu}f)\right\|_{q}^{q}&=\Ex{\Bigl|\Ex{(T_{\mu}f)(x_{k})\,\big|\,x_{0}}-(T_{\mu}f)(x_{0})\Bigr|^{q}}\nonumber\\
&\leq \Ex{\Bigl|(T_{\mu}(f))(x_{k})-(T_{\mu} f)(x_{0})\Bigr|^{q}},\label{eq:process_cumulate_difference}
\end{align}
where the last transition is due to the convexity of the function $t\mapsto |t|^{q}$. 
\item We know from~\Cref{prop:symmetry_D_mu_m} that the process is ``symmetric,'' in the sense that the joint distribution of $\left(x_{0},Z_{(1)},x_{1},\dots,Z_{(k)},x_{k}\right)$ is invariant under the reversal map
\[
\left(x_{0},Z_{(1)},x_{1},\dots,Z_{(k)},x_{k}\right)\mapsto \left(x_{k},Z_{(k)},x_{k-1},\dots,Z_{(1)},x_{0}\right).
\]
\item For each $j\in [k]$ and $i\in [m]$, let \[\Delta_{j,i}=f(z_{j,i})-(T_{\mu}f)(x_{j-1})\qquad\text{and}\qquad\Delta_{j,i}'=f(z_{j,i})-(T_{\mu}f)(x_{j}).\]
If we reveal the coordinates of each $Y_{(j)}$ one at a time, then $\Delta_{j,i}$ has conditional mean zero given the entire preceding history and the coordinates already revealed in the current step. Therefore,
    \[
        (\Delta_{1,1},\dots,\Delta_{1,m},
         \Delta_{2,1},\dots,\Delta_{2,m},
         \dots,\Delta_{k,1},\dots,\Delta_{k,m})
    \]
    is a martingale difference sequence in the sense of~\Cref{def:martingale_difference_sequence}. Since $x_{j-1}$ is uniform on $G$ and independent of $Y_{(j)}$, each $z_{j,i}$ is also uniformly distributed on $G$. Hence, by the triangle inequality and the fact that $T_\mu$ contracts the $L^q$-norm, we have
    \[
        \Ex{|\Delta_{j,i}|^q}^{1/q}
        \leq
        \|f\|_q+\|T_\mu f\|_q
        \leq
        2\|f\|_q.
    \]
We may then apply~\Cref{lem:martingale_n_step} and obtain
\begin{equation}\label{eq:martingale_forward_inequality}
\Ex{\left|\sum\nolimits_{j=1}^{k}\sum\nolimits_{i=1}^{m}\Delta_{j,i}\right|^{q}}\leq q^{q/2}\left(4km\|f\|_{q}^{2}\right)^{q/2}.
\end{equation}
By time-reversal symmetry, the sequence
\[
    (\Delta_{k,1}',\dots,\Delta_{k,m}',
     \Delta_{k-1,1}',\dots,\Delta_{k-1,m}',
     \dots,\Delta_{1,1}',\dots,\Delta_{1,m}')
\]
is also a martingale difference sequence. Thus,
\begin{equation}\label{eq:martingale_backward_inequality}
\Ex{\left|\sum\nolimits_{j=1}^{k}\sum\nolimits_{i=1}^{m}\Delta'_{j,i}\right|^{q}}\leq q^{q/2}\left(4km\|f\|_{q}^{2}\right)^{q/2}.
\end{equation}
\end{enumerate}

Combining~\Cref{eq:martingale_forward_inequality,eq:martingale_backward_inequality} (and using~\Cref{eq:q_norm_of_sum}), we have
\[
\Ex{\left|\sum\nolimits_{j=1}^{k}\sum\nolimits_{i=1}^{m}(\Delta_{j,i}-\Delta'_{j,i})\right|^{q}}^{1/q}\leq 4\sqrt{qkm}\,\|f\|_{q}.
\]
By the definition of the variables $\Delta_{j,i}$ and $\Delta_{j,i}'$, this can be rewritten as
\[
\Ex{\Bigl|(T_{\mu}(f))(x_{k})-(T_{\mu} f)(x_{0})\Bigr|^{q}}^{1/q}\leq 4\sqrt{\frac{qk}{m}}\,\|f\|_{q}.
\]
Combining with~\Cref{eq:process_cumulate_difference} then yields the desired conclusion.
\end{proof}

\subsection{Fourier-Analytic Perspectives}\label{subsec:Fourier_prelim}

As promised in~\Cref{subsec:bootstrapping}, we now return to the (vaguely stated) goal of ``replacing $T_{D(\mu,m)}$ by $T_{V}$ for some $V$,'' using a Fourier-analytic argument following~\cite{sanders2012bogolyubov}. The Fourier-analytic approach to almost-periodicity has a long history, including work of Bogolyubov~\cite{bogolyubov1939presque}, Bourgain~\cite{bourgain1990arithmetic}, and Chang~\cite{chang2002polynomial}; see also the discussion in~\cite[Section~9.2]{croot2010probabilistic}. In this subsection, we introduce this perspective and some related tools. In the next subsection, we combine them with~\Cref{lem:martingale_almost_period} to complete the argument.

We use the following notation for Fourier analysis on finite Abelian groups.

\begin{notation}
A character of $G$ is a group homomorphism $\chi:G\rightarrow \bC^{\times}$. Since $G$ is finite, we have $|\chi(x)|=1$ for all $x\in G$. We write $\widehat{G}$ for the group of characters of $G$, with pointwise multiplication as its group operation. For $f\in L^2(G)$ and $\chi\in\widehat G$, define
\[
    \widehat f(\chi)
    :=
    \Exu{x\in G}{f(x)\overline{\chi(x)}}.
\]
\end{notation}

Consider a function $f\in L^2(G)$, with Fourier expansion
\[
    f=\sum\nolimits_{\chi\in\widehat G}\widehat f(\chi)\chi.
\]
Applying a translation gives
\[
    T_{\{a\}}f
    =
    f(\cdot-a)
    =
    \sum\nolimits_{\chi\in\widehat G}
        \widehat f(\chi)\chi(\cdot-a)
    =
    \sum\nolimits_{\chi\in\widehat G}
        \left(\widehat f(\chi)\chi(a)^{-1}\right)\chi.
\]
Thus, $T_{\{a\}}$ multiplies the $\chi$-component of $f$ by $\chi(a)^{-1}$. Intuitively, if most of the Fourier mass of $f$ lies on characters $\chi$ for which $\chi(a)$ is close to $1$, then $T_{\{a\}}f$ should be close to $f$. This motivates the following standard definition.

\begin{definition}\label{def:bohr_set}
Given a character set $\Lambda\subseteq\widehat G$ and a parameter $\rho\in(0,2]$, the Bohr set with frequency set $\Lambda$ and radius $\rho$ is
\[
    B(\Lambda,\rho)
    :=
    \left\{
        x\in G\,\big|\,
        |\chi(x)-1|\leq\rho
        \text{ for every }\chi\in\Lambda
    \right\}.
\]
We say that $B(\Lambda,\rho)$ is a Bohr set with rank $|\Lambda|$ and radius $\rho$.
\end{definition}

Note that when $G=\bF_2^n$, Bohr sets are linear subspaces of $\bF_2^n$.

The following lemma makes the intuition preceding~\Cref{def:bohr_set} precise. It is known as a Bernstein-type inequality in harmonic analysis, and we include its proof in~\Cref{subsec:Bernstein} for the sake of completeness.

\begin{lemma}\label{lem:Bernstein}
Let $f\in L^2(G)$ have Fourier support contained in a character set $\Lambda\subseteq\widehat G$, meaning that $\widehat f(\chi)=0$ for every $\chi\notin\Lambda$. Then, for every $\rho\in(0,1)$ and every $a\in B(\Lambda,\rho)$, we have
\[
\left\|T_{\{a\}}f-f\right\|_{1}\leq 3\sqrt{2}\,\rho\|f\|_{1}.
\]
\end{lemma}

If $f$ has only a small number of nonzero Fourier coefficients, then~\Cref{lem:Bernstein} already gives a strong almost-periodicity result. In general, however, we cannot expect the Fourier support of $f$ to be small. What we can control is the set of Fourier coefficients with large magnitude. Suppose that $\|f\|_1=1$ and $\|f\|_{2}$ is not much larger. By Parseval's identity, the number of characters $\chi$ satisfying $\left|\widehat f(\chi)\right|\geq\eta$ is at most $\eta^{-2}\|f\|_2^2$. Chang's lemma~\cite{chang2002polynomial} provides a stronger structural description: these characters can be ``spanned'' by a collection of at most $C\eta^{-2}\log(2\|f\|_2)$ characters. The precise meaning of this span is defined below.

\begin{notation}
Given a function $f\in L^{2}(G)$ and a parameter $\eta\in(0,1]$, define
\[
    \operatorname{Spec}_\eta(f)
    :=
    \left\{
        \chi\in\widehat G\,\middle|\,
        \left|\widehat{f}(\chi)\right|\geq\eta
    \right\}.
\]
Given a character set $\Lambda\subseteq\widehat G$, write
\[
    \operatorname{Span}_{\pm}(\Lambda)
    :=
    \left\{
        \prod\nolimits_{\lambda\in\Lambda}\lambda^{\varepsilon_\lambda}\,\middle|\,
        \varepsilon_\lambda\in\{-1,0,1\}
        \text{ for every }\lambda\in\Lambda
    \right\}\subseteq \widehat{G}.
\]
\end{notation}

\begin{lemma}[Chang's lemma]\label{lem:chang_general}
Given a function $f\in L^{2}(G)$ with $\|f\|_{1}=1$ and a parameter $\eta\in(0,1]$, there exists a set $\Lambda\subseteq\operatorname{Spec}_\eta(f)$ such that
\begin{equation}\label{eq:Chang_lemma_requirement}
    |\Lambda|
    \leq
    32\eta^{-2}\log_{2}(2\|f\|_2)
    \qquad\text{and}\qquad
    \operatorname{Spec}_\eta(f)
    \subseteq
    \operatorname{Span}_{\pm}(\Lambda).
\end{equation}
\end{lemma}

We include a proof of this form of Chang's lemma in~\Cref{subsec:Chang}. The logarithmic dependence on $\|f\|_2$ is crucial to the argument of~\cite{sanders2012bogolyubov}: in (for instance)~\Cref{cor:dense_bogolyubov}, it allows us to obtain bounds on $|\Lambda|$ that are polynomial in $\log(1/\alpha)$  rather than polynomial in $1/\alpha$.

Now suppose we try to use~\Cref{lem:Bernstein,lem:chang_general} to prove almost-periodicity for a function of the form $g=f*f$, where $\|f\|_1=1$. The natural idea is to write $g=g_{\mathsf{in}}+g_{\mathsf{out}}$, where
\[
g_{\mathsf{in}}=\sum_{\chi\in \operatorname{Spec}_{\eta}(f)}\widehat{g}(\chi)\chi\qquad\text{and}\qquad g_{\mathsf{out}}=\sum_{\chi\in\widehat{G}\setminus \operatorname{Spec}_{\eta}(f)}\widehat{g}(\chi)\chi.
\]
We can use~\Cref{lem:chang_general} to find a small character set $\Lambda$ spanning $\operatorname{Spec}_\eta(f)$. If $a\in B(\Lambda,\rho)$, then every $\chi\in\operatorname{Spec}_\eta(f)$ satisfies
\[
    |\chi(a)-1|
    \leq
    \sum_{\lambda\in\Lambda}|\lambda(a)-1|
    \leq
    \rho|\Lambda|.
\]
Consequently, choosing $\rho$ sufficiently small allows us to apply~\Cref{lem:Bernstein} to $g_{\mathsf{in}}$, showing that every $a\in B(\Lambda,\rho)$ is an almost-period of $g_{\mathsf{in}}$. 

The remaining issue is to control $g_{\mathsf{out}}$, which by definition features only Fourier coefficients of small magnitude. The natural attempt via Parseval's identity gives
\[
    \|g_{\mathsf{out}}\|_1
    \leq
    \|g_{\mathsf{out}}\|_2
    =
    \left(
        \sum\nolimits_{\chi\notin\operatorname{Spec}_\eta(f)}
            \left|\widehat f(\chi)\right|^4
    \right)^{1/2}
    \leq
    \eta\|f\|_2.
\]
To obtain an error of order $\varepsilon$ from this estimate, however, we would have to take the threshold $\eta$ to be very small, on the order of $\varepsilon/\|f\|_2$. The factor $\eta^{-2}$ in~\Cref{eq:Chang_lemma_requirement} then leads to only a bound of order
\(
    \|f\|_2^{2+o(1)}
\) on $|\Lambda|$, instead of the desired logarithmic dependence on $\|f\|_{2}$. Eventually, in (for instance)~\Cref{cor:dense_bogolyubov}, this would translate to a bound on $|\Lambda|$ that is polynomial in $1/\alpha$.

In summary, the issue is that we do not have good ways to control the small Fourier coefficients, except for the crude upper bound given by Parseval's identity. \Cref{lem:martingale_almost_period} provides exactly the missing piece: the $k$-fold convolution (vaguely speaking) rapidly suppresses the small Fourier coefficients, so that we can apply~\Cref{lem:chang_general} with a larger threshold $\eta$. The detailed argument will be given in the next subsection.

\subsection{A Bohr Set of Almost-Periods}
\label{subsec:bohr_almost_periods}

Before combining~\Cref{lem:martingale_almost_period} with the Fourier tools, we first need to turn~\Cref{lem:martingale_almost_period} into an $L^{1}$ statement. This can be done using the duality between $L^{q}$-norms and $L^{q/(q-1)}$-norms.
\begin{lemma}\label{lem:martingale_almost_period_L1}
Let $m,k$ be positive integers, and let $\mu$ be a probability distribution on $G$. For any $K>1$ and any $f\in L^{2}(G)$ with $\|f\|_{2}\leq 2^{K}\|f\|_{1}$, we have
\[\left\|T_{\mu}f-T_{D(\mu,m)}^{k}(T_{\mu}f)\right\|_{1}\leq 16\sqrt{\frac{Kk}{m}}\,\|f\|_{1}.\]
\end{lemma}
\begin{proof}
For any probability distribution $\sigma$ on $G$, let $\sigma'$ be the distribution of $-a$ when $a$ is sampled from $\sigma$. It is easy to see from~\Cref{eq:def_of_T_mu} that $T_{\sigma'}$ is the adjoint of $T_{\sigma}$. It is also easy to see from~\Cref{def:D_mu_m} that $D(\sigma',m)=D(\sigma,m)'$. 

Now let $\nu=D(\mu,m)$. For any $g\in L^{2}(G)$ and any $q\geq 2$, we have
\[
\left|\inp{T_\mu f-T_\nu^k(T_\mu f)}{g}\right|=
\left|\inp{f}{T_{\mu'}g-T_{\nu'}^k(T_{\mu'}g)}\right|
\leq
4\sqrt{\frac{qk}{m}}\,
\|f\|_{q/(q-1)}\|g\|_q,
\]
where the last transition follows from H\"older's inequality and~\Cref{lem:martingale_almost_period}. Taking the supremem over all $g\in L^{2}(G)$ satisfying $\|g\|_{q}\leq 1$, we deduce that
\begin{equation}\label{eq:martingale_almost_period_dual}
\left\|T_{\mu}f-T_{D(\mu,m)}^{k}(T_{\mu}f)\right\|_{q/(q-1)}\leq 4\sqrt{\frac{qk}{m}}\,
\|f\|_{q/(q-1)}.
\end{equation}
By H\"older's inequality and the assumption on $f$, we have
\begin{equation}\label{eq:dual_interpolation}
\|f\|_{q/(q-1)}\leq
\|f\|_{1}^{1-2/q}\|f\|_{2}^{2/q}\leq 2^{2K/q}\|f\|_{1}.
\end{equation}
Combining~\Cref{eq:martingale_almost_period_dual,eq:dual_interpolation} and taking $q=2K$ (and using the fact that the $L^{1}$-norm of a function is at most its $L^{q/(q-1)}$-norm) yields the desired conclusion.
\end{proof}

We are now ready to apply the Fourier tools in~\Cref{subsec:Fourier_prelim} to prove our ultimate almost-periodicity lemma for $T_{\mu} f$. 

\begin{lemma}\label{lem:bohr_almost_period}
Let $K_1,K_2>1$ and $\varepsilon\in(0,1)$, and let $\mu$ be a probability distribution on $G$ satisfying $\|h_\mu\|_2\leq 2^{K_1}$. There exists a Bohr set $B(\Lambda,\rho)\subseteq G$ with
\[
    |\Lambda|
    \leq
    2^{22}\varepsilon^{-2}K_1K_2
    \bigl(K_2+\log_2(8/\varepsilon)\bigr)
    \qquad\text{and}\qquad
    \rho
    =
    \frac{\varepsilon}
         {50\max\{1,|\Lambda|\}},
\]
such that the following holds: for every $a\in B(\Lambda,\rho)$ and every $f\in L^2(G)$ satisfying $\|f\|_2\leq 2^{K_2}\|f\|_1$, we have
\begin{equation}\label{eq:bohr_almost_period}
    \left\|T_\mu f-T_{\{a\}}(T_\mu f)\right\|_1
    \leq
    \varepsilon\|f\|_1.
\end{equation}
\end{lemma}

\begin{proof}
Choose \(k=\left\lceil K_2+\log_2(8/\varepsilon)\right\rceil\) and \(m=\left\lceil2^{12}\varepsilon^{-2}K_{2}k\right\rceil.\) Plugging into~\Cref{lem:martingale_almost_period_L1}, we get
\begin{equation}\label{eq:bohr_smoothing_error}
\left\|T_{\mu}f-T_{D(\mu,m)}^{k}(T_{\mu}f)\right\|_{1}\leq \frac{\varepsilon}{4}\|f\|_{1}
\end{equation}
for any $f\in L^{2}(G)$ with $\|f\|_{2}\leq 2^{K_{2}}\|f\|_{1}$.

Let $\nu=D(\mu,m)$. Recall from~\Cref{lem:D_mu_m_density} that $\|h_\nu\|_2\leq2^{mK_1}$. Consider the character set $\Gamma=\operatorname{Spec}_{1/2}(h_{\nu})\subseteq \widehat{G}$. Applying~\Cref{lem:chang_general} with $\eta=1/2$, we obtain a set $\Lambda\subseteq\Gamma$ such that
\[
    \Gamma\subseteq\operatorname{Span}_{\pm}(\Lambda)\qquad\text{and}
    \qquad
    |\Lambda|\leq128(1+mK_1).
\]
Define $\rho$ as in the statement of this lemma. For any element $a\in B(\Lambda,\rho)$ and any character $\chi\in\Gamma$, the spanning property and the triangle inequality give
\(
    |1-\chi(a)|
    \leq |\Lambda|\rho
    \leq \varepsilon/50.
\)
Therefore, we have $B(\Lambda,\rho)\subseteq B(\Gamma,\varepsilon/50)$.

Now fix a function $f\in L^{2}(G)$ with $\|f\|_{2}\leq 2^{K_{2}}\|f\|_{1}$. We write $g=T_{\nu}^{k}(T_{\mu}f)$ and decompose it into $g=g_{\mathsf{in}}+g_{\mathsf{out}}$, where
\[
g_{\textsf{in}}=\sum\nolimits_{\chi\in \Gamma}\widehat{g}(\chi)\chi\qquad\text{and}\qquad g_{\textsf{out}}=\sum\nolimits_{\chi\in\widehat{G}\setminus \Gamma}\widehat{g}(\chi)\chi.
\]
By Parseval's identity,
\begin{align*}
    \left\|g_{\mathsf{out}}\right\|_2^2
    &=\sum\nolimits_{\chi\in \widehat{G}\setminus \Gamma}
        \left|\widehat{g}(\chi)\right|^2=\sum\nolimits_{\chi\in \widehat{G}\setminus \Gamma}
        \left|\widehat{h_\nu}(\chi)\right|^{2k}
        \left|\widehat{T_\mu f}(\chi)\right|^2\\
    &\leq 2^{-2k}\sum\nolimits_{\chi\in \widehat{G}}\left|\widehat{T_\mu f}(\chi)\right|^{2} \tag{by the definition of $\Gamma$}\\
    &=2^{-2k}
    \left\|T_\mu f\right\|_2^2\leq 2^{-2k}\|f\|_{2}^{2}\tag{$T_{\mu}$ contracts the $L^{2}$-norm}\\
    &\leq 2^{2K_{2}-2k}\|f\|_{1}^{2}\leq \left(\frac{\varepsilon}{8}\|f\|_{1}\right)^{2}. \tag{by the assumption on $f$ and the definition of $k$}
\end{align*}
Consequently, $\|g_{\mathsf{out}}\|_{1}\leq \|g_{\mathsf{out}}\|_{2}\leq \varepsilon\|f\|_{1}/8$, and hence
\[
\|g_{\mathsf{in}}\|_{1}\leq \|g\|_{1}+\|g_{\mathsf{out}}\|_{1}\leq \left\|T_{\nu}^{k}(T_{\mu}f)\right\|_{1}+\frac{\varepsilon}{8}\|f\|_{1}\leq 2\|f\|_{1}.
\]
For any $a\in B(\Lambda,\rho)$, since we also have $a\in B(\Gamma,\varepsilon/50)$, we may apply~\Cref{lem:Bernstein} and obtain
\[
\left\|T_{\{a\}}g_{\mathsf{in}}-g_{\mathsf{in}}\right\|_{1}\leq 3\sqrt{2}\cdot\frac{\varepsilon}{50}\cdot \|g_{\mathsf{in}}\|_{1}\leq \frac{\varepsilon}{4}\|f\|_{1}.
\]
Therefore, we have
\[
\left\|T_{\{a\}}g-g\right\|_{1}\leq \|T_{\{a\}}g_{\mathsf{out}}\|_{1}+\|g_{\mathsf{out}}\|_{1}+\left\|T_{\{a\}}g_{\mathsf{in}}-g_{\mathsf{in}}\right\|_{1}\leq \frac{\varepsilon}{8}\|f\|_{1}+\frac{\varepsilon}{8}\|f\|_{1}+\frac{\varepsilon}{4}\|f\|_{1}=\frac{\varepsilon}{2}\|f\|_{1}.
\]
Combining this with~\Cref{eq:bohr_smoothing_error} and translation invariance of the $L^1$-norm, we conclude that
\begin{align*}
    \left\|T_\mu f-T_{\{a\}}(T_\mu f)\right\|_1
    &\leq
    \|T_\mu f-g\|_1
    +\left\|g-T_{\{a\}}g\right\|_1
    +\left\|T_{\{a\}}(g-T_\mu f)\right\|_1\\
    &\leq
    \frac{\varepsilon}{4}\|f\|_1
    +\frac{\varepsilon}{2}\|f\|_1
    +\frac{\varepsilon}{4}\|f\|_1
    =\varepsilon\|f\|_1.\qedhere
\end{align*}
\end{proof}

\section{Applications}
\label{sec:applications}

In this section, we deduce~\Cref{thm:main_F2} and~\Cref{cor:dense_bogolyubov,cor:dense_difference_bohr} (easily) from~\Cref{lem:bohr_almost_period}.

\begin{proof}[Proof of~\Cref{thm:main_F2}]
Apply~\Cref{lem:bohr_almost_period} with the given $\mu$, $K_1$, $K_2$, and $\varepsilon$, and let $V=B(\Lambda,\rho)$ be the resulting Bohr set. Note that every character of $\bF_2^n$ takes values in $\{-1,1\}$, and the radius supplied by the lemma is less than $2$. Therefore,
\[
    V=\left\{x\in\bF_2^n\,\middle|\,\chi(x)=1\text{ for every }\chi\in\Lambda\right\}
\]
is a linear subspace of codimension at most $|\Lambda|$. The rank bound in~\Cref{lem:bohr_almost_period} then yields
\[
    \operatorname{codim}V
    \leq
    2^{22}\varepsilon^{-2}K_1K_2
    \bigl(K_2+\log_2(8/\varepsilon)\bigr)
    =O_\varepsilon(K_1K_2^2).
\]
For every $f$ satisfying $\|f\|_2\leq2^{K_2}\|f\|_1$, averaging~\Cref{eq:bohr_almost_period} over $a\in V$ gives
\[
    \|T_\mu f-T_V(T_\mu f)\|_1
    \leq
    \Exu{a\in V}{
        \left\|T_\mu f-T_{\{a\}}(T_\mu f)\right\|_1
    }\leq\varepsilon\|f\|_1.
\]
Since $V$ is chosen independently of $f$, this proves the theorem.
\end{proof}

For the two corollaries, we first record an elementary estimate for the size of a Bohr set in terms of $|\Lambda|$ and $\rho$; see, for example,~\cite[Lemma 4.20]{tao2006additive}.

\begin{lemma}\label{lem:bohr_size}
Let $G$ be a finite Abelian group. For every $\Lambda\subseteq\widehat G$ and $\rho\in(0,1]$, we have
\[
    |B(\Lambda,\rho)|
    \geq
    \left(\frac{\rho}{2\pi}\right)^{|\Lambda|}|G|.
\]
\end{lemma}

\begin{proof}
For each character $\chi\in\Lambda$, independently sample a random point $z_{\chi}$ uniformly on the unit circle $\{z\in\bC\mid|z|=1\}$. Consider the set
\(
S=\left\{x\in G\,\big|\,\left|\chi(x)-z_{\chi}\right|\leq \rho/2\text{ for all }\chi\in\Lambda\right\}\subseteq G
\). It is easy to see that 
\[\Exu{z_{\chi}}{|S|}\geq \left(\frac{\rho}{2\pi}\right)^{|\Lambda|}|G|,\]
so there exists a choice of the points $z_{\chi}$ for which $|S|\geq (\rho/(2\pi))^{|\Lambda|}|G|$. Fixing such a set $S$, if we pick any $x_{0}\in S$, we have $|\chi(x-x_{0})-1|=|\chi(x)-\chi(x_{0})|\leq \rho$ for any $\chi\in\Lambda$ and $x\in S$. Therefore, the translation $S-x_{0}$ is contained in $B(\Lambda,\rho)$, and hence $|B(\Lambda,\rho)|\geq (\rho/(2\pi))^{|\Lambda|}|G|$.
\end{proof}

\begin{proof}[Proof of~\Cref{cor:dense_bogolyubov}]
Let $g=h_{A}*h_{A}=T_{A}h_{A}$, and recall that $L=\log_2(2/\alpha)$. Apply~\Cref{lem:bohr_almost_period} with $\varepsilon=1/2$ and
\begin{equation}\label{eq:set_parameters}
\mu\text{ uniform on }A,\qquad f=h_{A},\qquad K_{1}=K_{2}=2+\log_{2}(1/\alpha)/2.
\end{equation}
This gives a Bohr set $B=B(\Lambda,\rho)$ with $|\Lambda|\leq CL^{3}$ and $\rho\geq C^{-1}L^{-3}$, where $C>0$ is an absolute constant, such that for each $b\in B$ we have $\left\|g-T_{\{b\}}g\right\|\leq 1/2$. In particular, the supports of the probability density functions $g$ and $T_{\{b\}}g$ cannot be disjoint. Thus, $2A\cap (2A+b)\neq \emptyset$, which implies $b\in 2A-2A$. We therefore conclude that $B\subseteq 2A-2A$.

The bounds on $|\Lambda|$ and $\rho$, together with~\Cref{lem:bohr_size}, imply that
\[
    |B|
    \geq
    \exp\!\left(-|\Lambda|\log\frac{2\pi}{\rho}\right)|G|
    \geq
    \exp\bigl(-C'L^3\log(2L)\bigr)|G|
\]
for some absolute constant $C'>0$. If $G=\bF_2^n$, the same argument as in the proof of~\Cref{thm:main_F2} shows that $B$ is a linear subspace of codimension at most $CL^3$.
\end{proof}

\begin{proof}[Proof of~\Cref{cor:dense_difference_bohr}]
Let $g=h_{A}*h_{A}=T_{A}h_{A}$. Apply~\Cref{lem:bohr_almost_period} with the given constant $\varepsilon$ and the other parameters set as in~\Cref{eq:set_parameters}. 
This gives a Bohr set $B=B(\Lambda,\rho)$ with $|\Lambda|\leq C_{\varepsilon}L^{3}$ and $\rho\geq C_{\varepsilon}^{-1}L^{-3}$, where $C_{\varepsilon}>0$ depends only on $\varepsilon$, such that for each $b\in B$ we have $\left\|g-T_{\{b\}}g\right\|\leq \varepsilon$. Averaging over $b\in B$ and using the convexity of the $L^1$-norm, we obtain
\[
    \|g-T_Bg\|_1\leq\varepsilon.
\]
Since $g$ is a probability density function supported on $2A$, this means that an element sampled according to the density function $T_Bg$ lies in $2A$ with probability at least $1-\varepsilon$. Equivalently,
\[
\Pru{a_{1},a_{2}\in A,\,b\in B}{a_{1}+a_{2}+b\in 2A}\geq 1-\varepsilon
.\]
where $a_1,a_2$, and $b$ are sampled independently and uniformly from their respective sets. Thus there exists a realization of $a_{1}+a_{2}$ for which
\[
\frac{|(a_{1}+a_{2}+B)\cap 2A|}{|B|}\geq 1-\varepsilon,
\]
Taking $x=a_1+a_2\in2A$ gives the desired inequality $|(x+B)\cap 2A|\geq (1-\varepsilon)|B|$. The rest follow as in the proof of~\Cref{cor:dense_bogolyubov}.
\end{proof}

\addcontentsline{toc}{section}{References}
\bibliographystyle{alpha}
\bibliography{reference}

\appendix
\section{Probability Tools}\label{sec:sum_fluctuations}

In this appendix, we present some basic probability inequalities used in the body of the paper.

The following (well-known) inequality is used in~\Cref{subsec:empirical_approx}.
\begin{proposition}[Marcinkiewicz--Zygmund inequality]
\label{prop:Marcinkiewicz-Zygmund}
Let $q\geq2$ be an even integer, and let $X_1,\dots,X_n$ be independent complex-valued random variables satisfying $\Ex{|X_i|^q}<\infty$ for every $i\in[n]$. Then
\[
    \Ex{
        \left|\sum\nolimits_{i=1}^n\bigl(X_i-\Ex{X_i}\bigr)\right|^q
    }
    \leq
    (2q)^{q/2}\cdot
    \Ex{
        \left(\sum\nolimits_{i=1}^n|X_i|^2\right)^{q/2}
    }.
\]
\end{proposition}

\begin{proof}
If the tuple of signs $\sigma=(\sigma_{1},\dots,\sigma_{n})$ is uniformly random in $\{-1,1\}^{n}$, the inequality
\[
    \Exu{\sigma}{
        \left|\sum\nolimits_{i=1}^n\sigma_i z_i\right|^{q}
    }
    \leq
    (q/2)^{q/2}\left(\sum\nolimits_{i=1}^n|z_i|^2\right)^{q/2}
    \qquad\text{for all }z_1,\dots,z_n\in\bC.
\]
can be proved by expanding the left-hand side, taking absolute values of the surviving terms, and then comparing coefficients with those in the expansion of the right-hand side.

Now let $X_1',\dots,X_n'$ be independent copies of $X_1,\dots,X_n$, independent also of $\sigma$. Conditional Jensen's inequality and the symmetry of each difference $X_i-X_i'$ give
\begin{align*}
    \Ex{
        \left|\sum\nolimits_{i=1}^n\bigl(X_i-\Ex{X_i}\bigr)\right|^{q}
    }
    &\leq
    \Ex{
        \left|\sum\nolimits_{i=1}^n\bigl(X_i-X_{i}'\bigr)\right|^{q}
    }=
    \Ex{
        \left|\sum\nolimits_{i=1}^n\sigma_{i}\bigl(X_i-X_{i}'\bigr)\right|^{q}
    }\\
    &\leq
    2^{q}\cdot
    \Ex{
        \left|\sum\nolimits_{i=1}^n\sigma_iX_i\right|^{q}
    }\leq
    (2q)^{q/2}\cdot 
    \Ex{
        \left(\sum\nolimits_{i=1}^n|X_i|^2\right)^{q/2}
    }.
\end{align*}
The second-to-last inequality uses
\begin{equation}\label{eq:q_norm_of_sum}
    |u-v|^{q}
    \leq
    2^{q-1}\bigl(|u|^{q}+|v|^{q}\bigr)\qquad\text{for all }u,v\in \bC
\end{equation}
and the fact that $(\sigma_{1}X_1,\dots,\sigma_{n}X_n)$ and $(\sigma_{1}X_1',\dots,\sigma_{n}X_n')$ have the same distribution.
\end{proof}

\subsection{Martingale Increments}\label{subsec:martingale_increment}

In~\Cref{subsec:martingale}, we come across sums of mean-zero random variables that need not be fully independent, but still form the ``difference sequence'' of a martingale. The precise definition is given below.

\begin{definition}\label{def:martingale_difference_sequence}
Let $X_1,\dots,X_n$ be jointly distributed random vectors in $\bR^{m}$ satisfying $\Ex{|X_i|}<\infty$ for every $i\in[n]$. We say that $(X_1,\dots,X_n)$ form a martingale difference sequence if
\[
    \Ex{X_i\mid X_1,\dots,X_{i-1}}=0
    \qquad\text{almost surely for every }i\in[n].
\]
For $i=1$, this condition means $\Ex{X_1}=0$.
\end{definition}

Although we are concerned with complex-valued random variables in~\Cref{subsec:martingale}, among the complex number operations we use only addition, real scalar multiplication, and the Euclidean norm. We can therefore identify the complex plane $\bC$ with the real Euclidean space $\bR^2$, and we can in fact work more generally with random vectors in $\bR^m$.

\begin{notation}
For a random vector $X$ and a real number $q\geq 1$, write
\[
    \|X\|_q:=\left(\Ex{|X|^q}\right)^{1/q}.
\]
\end{notation}

To analyze the $p$-th power of the norm of a sum of two complex numbers, we use the following elementary lemma.

\begin{lemma}\label{lem:vector_Hessian}
Let $a,b\in \bR^{m}$ be vectors and let $q\geq 2$ be a real number. Then the function $f:\bR\rightarrow \bR$ defined by $f(t)=|a+tb|^{q}$ is twice continuously differentiable, and we have
\[
f'(0)=q|a|^{q-2}\langle a,b\rangle,\qquad
0\leq f''(t)\leq q(q-1)|a+tb|^{q-2}|b|^{2}\qquad\text{for all }t\in \bR,
\]
where for $q=2$ the factor $|a|^{0}$ and $|a+tb|^{0}$ are understood to be constantly 1.
\end{lemma}
\begin{proof}
Although these conclusions follow from direct calculations, we provide some geometric intuition here. Consider the map $F:\bR^{m}\rightarrow\bR$ given by $F(z)=|z|^q$. Its gradient at $z\in \bR^{m}$ is $q|z|^{q-2}z$, pointing in the radial direction. At any point $(x,0,\dots,0)$ with $x\in\bR\setminus\{0\}$, the Hessian of $F$ is
\[
\begin{bmatrix}
q(q-1)|x|^{q-2} & 0\\
0 & q|x|^{q-2}\cdot I_{m-1}
\end{bmatrix}.
\]
By rotational invariance, it is then easy to see that $F$ is twice continuously differentiable (and convex) as a real map. Furthermore, at every nonzero point $z$, the second directional derivative among unit directions is maximized in the radial direction, where it equals $q(q-1)|z|^{q-2}$. Applying this to the line $t\mapsto a+tb$, whose velocity has length $|b|$, gives the desired inequality.
\end{proof}

To bound the $q$-norm of a sum of martingale differences, we first prove the following estimate for two increments.

\begin{lemma}\label{lem:martingale_one_step}
Let $q\geq 2$ be a real number, and let $(X_{1},X_{2})$ be a martingale difference sequence such that $\|X_{1}\|_{q},\|X_{2}\|_{q}<\infty$. We have
\[
\|X_{1}+X_{2}\|_{q}^{2}\leq\|X_{1}\|_{q}^{2}+(q-1)\|X_{2}\|_{q}^{2} .
\]
\end{lemma}
\begin{proof}
If $\|X_{1}\|_{q}=0$ this is trivially true. We henceforth assume $\|X_{1}\|_{q}>0$. Define functions $f,g:\bR\rightarrow \bR$ by letting
\[
f(t)=\|X_{1}+tX_{2}\|_{q}^{q}\qquad\text{and}\qquad g(t)=\|X_{1}+tX_{2}\|_{q}^{2}=f(t)^{2/q}.
\]
By the dominated convergence theorem, we can differentiate under the expectation in the following calculation:
\[
f'(t)=\frac{\d}{\d t}\Ex{|X_{1}+tX_{2}|^{q}}=\Ex{\frac{\d}{\d t}|X_{1}+tX_{2}|^{q}}
\]
Since $\Ex{X_{2}|X_{1}}=0$, it follows from~\Cref{lem:vector_Hessian} that
\[
\Ex{\frac{\d}{\d t}\Big|_{t=0}|X_{1}+tX_{2}|^{q}\,\middle|\,X_{1}}=0,
\]
and hence $f'(0)=0$. Differentiating a second time and applying~\Cref{lem:vector_Hessian}, we get
\[
0\leq f''(t)\leq q(q-1)\Ex{|X_{1}+tX_{2}|^{q-2}|X_{2}|^{2}}\leq q(q-1)f(t)^{(q-2)/q}\cdot\|X_{2}\|_{q}^{2},
\]
where we used H\"older's inequality in the last transition. Since $f''\geq0$ and $f'(0)=0$, we have $f(t)\geq f(0)>0$ for every $t\in\bR$. We then have $g'(0)=\frac{2}{q}f(0)^{(2-q)/q}f'(0)=0$, and for all $t\in \bR$,
\[
g(t)''=\frac{2}{q}f(t)^{(2-q)/q}f''(t)-\frac{2(q-2)}{q^{2}}f(t)^{(2-2q)/q}f'(t)^{2}\leq \frac{2}{q}f(t)^{(2-q)/q}f''(t)\leq 2(q-1)\|X_{2}\|_{q}^{2}.
\]
Therefore, 
\[g(1)-g(0) =
    \int_0^1(1-t)g''(t)\,\d t
    \leq (q-1)\|X_{2}\|_{q}^{2},\]
as desired.
\end{proof}

We can now iterate~\Cref{lem:martingale_one_step} to handle general martingale difference sequences. Compared with~\Cref{prop:Marcinkiewicz-Zygmund} in the fully-independent case, the following estimate has a better-optimized leading constant, but replaces $\Ex{(\sum_i|X_i|^2)^{q/2}}$ by the generally larger quantity $(\sum_i\|X_i\|_q^2)^{q/2}$.

\begin{lemma}\label{lem:martingale_n_step}
Let $q\geq 2$ be a real number, and let $(X_{1},X_{2},\dots,X_{n})$ be a martingale difference sequence such that $\|X_{i}\|_{q}< \infty$ for all $i\in [n]$. Then
\[
\Ex{\left|\sum\nolimits_{i=1}^{n}X_{i}\right|^{q}}\leq (q-1)^{q/2}\cdot\left(\sum\nolimits_{i=1}^{n}\Ex{|X_{i}|^{q}}^{2/q}\right)^{q/2}.
\]
\end{lemma}
\begin{proof}
The conclusion is the same as $\|X_{1}+\dots+X_{n}\|_{q}^{2}\leq (q-1)\left(\|X_{1}\|_{q}^{2}+\dots+\|X_{n}\|_{q}^{2}\right)$, which easily follows from~\Cref{lem:martingale_one_step} by induction on $n$.
\end{proof}

\section{Fourier Analysis Tools}
\label{sec:Fourier_tools}

\subsection{Bernstein's Inequality}\label{subsec:Bernstein}

To prove~\Cref{lem:Bernstein}, we first consider the case $G=\bZ/N\bZ$. The following proof is due to~\cite{bloom1974bernstein}.

\begin{lemma}\label{lem:Bernstein_cyclic}
Let $N\geq 2$ be a positive integer and let $G=\bZ/N\bZ$. Suppose $f\in L^{2}(G)$ satisfies $\widehat{f}(\chi)=0$ for all characters $\chi\in \widehat{G}$ such that $|\chi(1)-1|> \rho$, where $\rho\in (0,1)$ is a real number. Then we have
\[
\left\|T_{\{1\}}f-f\right\|_{1}\leq 3\sqrt{2}\,\rho\|f\|_{1}.
\]
\end{lemma}
\begin{proof}
In this proof, we let $i$ denote the imaginary unit. We list the characters in $\widehat{G}$ as $\{\chi_{j}\}_{j\in \bZ/N\bZ}$, where $\chi_{j}(x)=e^{2\pi ijx/N}$ for any $j,x\in \bZ/N\bZ$. Let $E=\left\{j\in \bZ/N\bZ\,\big|\,|\chi_{j}(1)-1|\leq \rho\right\}$. By assumption, the Fourier support of $f$ is contained in $\{\chi_j:j\in E\}$. Since $\rho<1$, we know that $E=\{-r,-r+1,\dots,r\}$ for some nonnegative integer $r\leq N/6$.

\paragraph{Step 1: the idea of convolution kernels.} We know that the Fourier expansion
\[
f=\sum\nolimits_{j\in E}\widehat{f}(\chi_{j})\chi_{j},
\]
decomposes $f$ into summands that each change by a multiplicative factor close to 1 under the action of $T_{\{1\}}$. It is then easy to see that
\[
\left\|T_{\{1\}}f-f\right\|_{1}\leq \rho\sum\nolimits_{j\in E}\left|\widehat{f}(\chi_{j})\right|.
\]
However, the Fourier $\ell^{1}$-norm on the right hand side need not be bounded by $\|f\|_{1}$, so this decomposition does not lead to the desired conclusion. The idea is to still decompose $f$ into pieces whose changes under $T_{\{1\}}$ can be controlled, but with the sum of the absolute values of the coefficients bounded by $\|f\|_1$.

Fortunately, the assumption on the Fourier support of $f$ also gives us room to choose such a decomposition. Indeed, for any function $g\in L^{2}(G)$ such that $\widehat{g}(\chi_{j})=1$ for all $j\in E$, we have $f*g=f$. In other words, we can write
\[
f=\sum_{x\in G}\frac{f(x)}{|G|}g(\cdot-x),
\]
decomposing $f$ into a linear combination of translates of $g$, with the sum of the absolute values of the coefficients equal to $\|f\|_1$. Therefore, it suffices to choose such a function $g$ with $\|T_{\{1\}}g-g\|_{1}\leq 3\sqrt{2}\rho$.

\paragraph{Step 2: finding a convolution kernel.} The idea is to find a function $g$ of the form $g=u_{1}u_{2}$, where $u_{1}$ and $u_{2}$ are two functions such that the $L^{2}$-norm of $\left(T_{\{1\}}u_{1}-u_{1}\right)$ and $\left(T_{\{1\}}u_{2}-u_{2}\right)$ can be bounded, so that we can then bound the $L^{1}$-norm of $\left(T_{\{1\}}g-g\right)$ by Cauchy-Schwarz. We let
\[
u_{1}=|E|^{-1/2}\sum\nolimits_{j\in E}\chi_{j}\qquad\text{and}\qquad u_{2}=|E|^{-1/2}\sum\nolimits_{j\in F}\chi_{j},
\]
where $F=\{-2r,-2r+1\dots,2r\}\subseteq \bZ/N\bZ$. Then, letting $g=u_{1}u_{2}$, we easily see that $\widehat{g}(\chi_{j})=1$ for all $j\in E$. Furthermore, 
\begin{align*}
\left\|T_{\{1\}}u_{1}-u_{1}\right\|_{2}^{2}&=\frac{1}{|E|}\sum_{j\in E}|\chi_{j}(1)-1|^{2}\leq \rho^{2},\qquad\text{and}\\
\left\|T_{\{1\}}u_{2}-u_{2}\right\|_{2}^{2}&=\frac{1}{|E|}\sum_{j\in F}|\chi_{j}(1)-1|^{2}\leq 8\rho^{2}.
\end{align*}
Therefore, we have
\begin{align*}
\left\|T_{\{1\}}g-g\right\|_{1}&\leq \left\|\left(T_{\{1\}}u_{1}-u_{1}\right)u_{2}\right\|_{1}+\left\|\left(T_{\{1\}}u_{1}\right)\left(T_{\{1\}}u_{2}-u_{2}\right)\right\|_{1}\\
&\leq \left\|T_{\{1\}}u_{1}-u_{1}\right\|_{2}\|u_{2}\|_{2}+\left\|T_{\{1\}}u_{1}\right\|_{2}\left\|T_{\{1\}}u_{2}-u_{2}\right\|_{2}\\
&\leq \rho\cdot \sqrt{2}+1\cdot 2\sqrt{2}\rho=3\sqrt{2}
\rho.\end{align*}
It then follows that
\[
\left\|T_{\{1\}}f-f\right\|_{1}=\left\|T_{\{1\}}(f*g)-(f*g)\right\|_{1}=\left\|f*\left(T_{\{1\}}g-g\right)\right\|_{1}\leq \|f\|_{1}\left\|T_{\{1\}}g-g\right\|_{1}\leq 3\sqrt{2}\,\rho\|f\|_{1}.\qedhere
\]
\end{proof}

We now deduce the result for general finite Abelian groups by applying the cyclic case on each coset of the subgroup generated by $a$.

\begin{proof}[Proof of~\Cref{lem:Bernstein}]
We can assume $a$ is not the identity element of the group $G$, and let $N\geq 2$ be the order of $a$. Partition $G$ into cosets of the subgroup $A=\{0,a,\dots,(N-1)a\}$. For each coset $x+A=\{x,x+a,\dots,x+(N-1)a\}$, consider the affine map $\iota:\bZ/N\bZ\rightarrow (x+A)$ defined by $\iota(j)=x+ja$. For every $\chi\in\widehat G$, we have
\[
    (\chi\circ\iota)(j)=\chi(x)\chi(a)^j.
\]
Thus, $\chi\circ\iota$ is a character of $\bZ/N\bZ$ multiplied by the constant $\chi(x)$, whose absolute value is $1$. Since the characters in the Fourier support of $f$ satisfy $|\chi(a)-1|\leq\rho$, the function $f\circ\iota$ satisfies the hypotheses of~\Cref{lem:Bernstein_cyclic}. Applying that lemma on each coset and then averaging over the cosets yields the desired conclusion.
\end{proof}

\subsection{Chang's Lemma}\label{subsec:Chang}

We remark that Chang's lemma is closely related to the ``level-1 inequality'' in the analysis of Boolean functions, and the following proof is similar to the standard Fourier-analytic proof of hypercontractivity on the Boolean hypercube.

\begin{proof}[Proof of~\Cref{lem:chang_general}] We divide the proof into the following steps.

\paragraph{Step 1: finding $\Lambda$.} A character set $\Lambda\subseteq \widehat{G}$ is said to be dissociated if for any signs $\varepsilon_{\lambda}\in \{-1,0,1\}$ that are not all zeros, the product $\prod_{\lambda\in \Lambda}\lambda^{\varepsilon_{\lambda}}$ does not equal the trivial character $1$. In particular, a dissociated character set cannot contain the trivial character. We claim that any dissociated characeter set $\Lambda$ that is inclusion-maximal among subsets of $\operatorname{Spec}_{\eta}(f)$ satisfies both conditions in~\Cref{eq:Chang_lemma_requirement}.

Indeed, maximality easily implies the required spanning property. Consider an arbitrary character $\chi\in \operatorname{Spec}_{\eta}(f)\setminus\Lambda$. Since $\Lambda\cup\{\chi\}$ is not dissociated, there is a nontrivial relation $\chi^{\varepsilon_{0}}\prod_{\lambda\in\Lambda}\lambda^{\varepsilon_{\lambda}}=1$, where $\varepsilon_{0},\varepsilon_{\lambda}\in\{-1,0,1\}$. We must have $\varepsilon_{0}\in\{-1,1\}$, because otherwise this would be a nontrivial relation among the elements of $\Lambda$. Therefore $\chi=\prod_{\lambda\in\Lambda}\lambda^{-\varepsilon_{0}\varepsilon_{\lambda}}\in \operatorname{Span}_{\pm}(\Lambda)$.

\paragraph{Step 2: H\"older's inequality.} In the remainder of the proof, we show the required upper bound on $|\Lambda|$ for an arbitrary dissociated character set $\Lambda\subseteq \widehat{G}$.

Consider a function $g\in L^{2}(G)$ defined by
\(
g(x)=\sum_{\lambda\in \Lambda}a_{\lambda}\lambda(x),
\)
where $a_{\lambda}\in \bC$ is such that $|a_{\lambda}|=1$ and $\widehat{f}(\lambda)\overline{a_{\lambda}}=\left|\widehat{f}(\lambda)\right|$, for each $\lambda\in \Lambda$. By Parseval's identity, we have
\begin{equation}\label{eq:Chang_Parseval}
\langle f,g\rangle=\sum\nolimits_{\lambda\in\Lambda}\widehat{f}(\lambda)\overline{a_{\lambda}}=\sum\nolimits_{\lambda\in\Lambda}\left|\widehat{f}(\lambda)\right|\geq \eta |\Lambda|.
\end{equation}
On the other hand, for any a real number $q\geq 2$, applying H\"older's inequality twice gives
\begin{equation}\label{eq:Chang_Holder}
|\langle f,g\rangle|\leq \|f\|_{q/(q-1)}\|g\|_{q}\leq \|f\|_{1}^{1-2/q}\|f\|_{2}^{2/q}\|g\|_{q}=\|f\|_{2}^{2/q}\|g\|_{q}.
\end{equation}
Combining~\Cref{eq:Chang_Parseval,eq:Chang_Holder}, we see that the natural next goal is to show a good upper bound on $\|g\|_q$ for large values of $q$.

\paragraph{Step 3: bounding $\|g\|_{q}$.} Thinking of $x$ as uniformly random on $G$, we view $g(x)$ as a sum of the random variables $a_{\lambda}\lambda(x)$, for $\lambda\in \Lambda$. Each of these random variables (and any square-free product of them) has mean 0, so it's natural to wonder whether the tools in~\Cref{sec:sum_fluctuations} can be used to bound the $q$-norm of their sum. However, these random variables need not be independent (or arranged as a martingale difference sequence). A standard trick in such context is to ``compare'' these random variables with independent signs (cf. the proof of~\Cref{prop:Marcinkiewicz-Zygmund}), in order to capture the cancellations among them. 

Let 
\[u_{\lambda}(x)=\operatorname{Re} \left(a_{\lambda}\lambda(x)\right)=\frac{a_{\lambda}\lambda(x)+\overline{a_{\lambda}}\lambda(x)^{-1}}{2}\in [-1,1]\qquad\text{for each }\lambda\in \Lambda\text{ and }x\in G.\] 
By expanding the products and using the condition that $\Lambda$ is dissociated, we have
\begin{equation}\label{eq:dissociated_u}
\Exs{x\in G}{\prod\nolimits_{\lambda\in\Gamma}u_{\lambda}(x)}=0\qquad\text{for all nonempty }\Gamma\subseteq \Lambda.
\end{equation}
Now consider the following process:
\begin{enumerate}
\item Sample a $x\in G$ uniformly at random.
\item Conditioned on $x$, independently sample signs $\sigma_{\lambda}\in\{-1,1\}$ according to $\Ex{\sigma_{\lambda}\,\middle|\,x}=u_{\lambda}(x)$.
\end{enumerate}
For any $s\in \{-1,1\}^{\Lambda}$, we have
\[
\Pr{\sigma_{\lambda}=s_{\lambda}\text{ for all }\lambda\in\Lambda}=\Exu{x\in G}{\prod\nolimits_{\lambda\in \Lambda}\frac{1+s_{\lambda}u_{\lambda}(x)}{2}}=2^{-|\Lambda|},
\]
where the last transition is due to~\Cref{eq:dissociated_u}. Therefore, after averaging over $x\in G$, the signs $\sigma_{\lambda}$ are independent and uniformly distributed. Since
\[
\operatorname{Re}(g(x))=\sum\nolimits_{\lambda\in\Lambda}u_{\lambda}(x)=\Ex{\sum\nolimits_{\lambda\in\Lambda}\sigma_{\lambda}(x)\,\middle|\,x},
\]
conditional Jensen’s inequality gives
\[
\|\operatorname{Re}(g)\|_{q}\leq \Ex{\left|\sum\nolimits_{\lambda\in \Lambda}\sigma_{\lambda}\right|^{q}}^{1/q}\leq\sqrt{(q-1)|\Lambda|}, 
\]
where the last transition is due to~\Cref{lem:martingale_n_step}. Similarly, we have the same upper bound on $\|\operatorname{Im}(g)\|_{q}$, and combining with~\Cref{eq:Chang_Parseval,eq:Chang_Holder} yields
\[
|\Lambda|\leq 2\eta^{-1}\|f\|_{2}^{2/q}\sqrt{(q-1)|\Lambda|}.
\]
Choosing $q=2\log_{2}(2\|f\|_{2})\geq 2$ now yields
\(
|\Lambda|\leq 32\eta^{-2}\log_{2}(2\|f\|_{2})
\), as desired.
\end{proof}
\end{document}